
\def\p {\partial}
\def \D {\Delta_{d{\bf v}}}
\def \Ds  {\Delta^{\#}}
\def\t{\tilde}
\def\s {\sigma}
\def\L {\Lambda}

\def\a{\alpha}
\def\O{\Omega}
\def\d{\delta}
\def\dv  {{d{\bf{v}}}}
\def\A {{\cal K}}
\def \PTM {\Pi T^*M}
\def\F{{\cal F}}
\def \AA {{\cal A}}
\def\O{\Omega}
\def\E {{\cal E}}
\def\G{{\cal G}}
\def\AA{{\cal A}}
\def\B{{\cal B}}

\def\S{{\cal B}}
\def\Ocal {{A}} 

\def\identifyl{\varphi_{_L}}

\def\ltimes {{{\vrule height1.7mm depth0mm width 0.5pt}}\!\times}

\def\finish {{\vrule height2mm depth0mm width 4pt}}
\def\superspace{supermanifold $\,$}


\def\BarScha{1}
\def\BarSchb{2}
\def\tyutina {3}
\def\tyutinb {4}
\def \BVa {5}
\def \BVb {6}
\def\Berezin {7}
\def\Bernsteina {8}
\def\Bernsteinb {9}
\def\buttin{10}
\def \Gayduk {11}
\def\Gui  {12}
\def\JMPa {13}
\def\CMP {14}
\def\But {15}
\def\MPL{16}
\def\JMPd {17}
\def\Poin {18}
\def\leitb {19}
\def\leites{20}
\def\Shander{21}
\def\SchNucl {22}
\def\SchCMPb {23}
\def\SchCMPa{24}
\def\Vora   {25}
\def\Vor  {26}


\def\statement {1.3}
\def\statementa {1.3a}
\def\statementb {1.3b}
\def\statementc {1.3c}
\def\cartan {1.4}


\def\bracket {2.1}
\def\hamvectorfield {2.2}
\def\jacoby {2.3}
\def\darbouxtheorem {2.4}
\def\berdef{2.6}
\def\deltadefgen{2.7}
 \def\deltadef {2.8}
\def\deltaprop {2.10}
\def\deltatransform {2.11}
\def\deltadensdef {2.12}
\def\correctness {2.13}
\def\cantransform {2.14}
\def\hamiltoniandeformation {2.15}
\def \denstransform {2.16}
\def\denscalculus {2.17}
\def\deltadenspropa {2.18}
\def\deltadenspropb {2.19}
\def\batsemidensity {2.20}
\def\bvsdensdef{2.21}


\def\diffdenscorr {3.4}
\def\bercanonical{3.5}
 \def\innerdelta {3.6}
\def\divdelta {3.7}
\def\innerdeltab {3.8}
\def\dopodin{3.9}
\def\firstoperation {3.10}
\def\secondoperation {3.12}

\def\identsympl {4.1}
\def\sdensforms {4.2}
\def\adjustedatlas {4.3}

\def\adjincoord {4.4}
\def\implies {4.5}
\def\hamdva{4.6}

\def\projectionmap{4.8}
\def\newdva{4.9}
\def\projectionchanging {4.10}
\def\projectionmaphom{4.11}

\def\spaceV{4.12}
\def\mapA {4.13}
\def\propertyofmapA{4.14}
\def\hamtimedepend {4.15}
\def\volumemap{4.16}
\def\changingdva {4.17}
\def\spaceVcoh{4.18}
\def\mapcoh {4.19}
\def\bijectivity {4.20}
\def\cohomologymap {4.21}
\def\cohmapinvariant {4.22}
\def\svjazjform {4.23}
\def\anothersurface{4.24}
\def\svjazjformcom {4.25}
\def\twosurfacefunction {4.26}
\def\svjazjforma {4.27}

\def\bvodin{4.28}
\def\bvsdensnormalform {4.29}
\def\bvdva{4.30}


\def\odin {5.1}
\def\dva{5.2}
\def\tri{5.3}
 \def\four{5.4}
\def\five{5.5}

 \def\theorem {5.7}
\def\flag{5.8}
\def\newdensities {5.9}
\magnification=1200 \baselineskip 14pt

\def\s {\sigma}
\def\L {\Lambda}

\def\a{\alpha}
\def\O{\Omega}
\def\d{\delta}
\def\dv  {{d{\bf{v}}}}
\def\A {{\cal K}}
\def \K {{\cal K}}
\def \PTM {\Pi T^*M}
\def\S{{\cal B}}
\def\F{{\cal F}}
\def \AA {{\cal A}}
\def\O{\Omega}

\def\identifyl{\varphi_{_L}}

\def\finish {{\vrule height2mm depth0mm width 4pt}}

\def\statement {1.3}
\def\statementa {1.3a}
\def\statementb {1.3b}
\def\statementc {1.3c}

                 \medskip

     \centerline {\bf Semidensities
            on Odd Symplectic Supermanifolds}

                      \medskip

                   \centerline {\bf Hovhannes (O.M.) Khudaverdian}
               \medskip

\centerline {\it Department of Mathematics}
 \centerline {\it UMIST, Manchester M60 1QD, UK}
       \centerline {\it and}
\centerline {\it Laboratory of Computing Technique and
        Automation}
\centerline   {\it Joint Institute for Nuclear Research,}

\centerline   {\it Dubna 141980, Russia}

\smallskip

 \centerline {\it on leave of absence of G.S. Sahakyan
        Department of Theoretical
       Physics }
  \centerline{\it of Yerevan State University,
    A. Manoukian St., 375049  Yerevan, Armenia}

\noindent e-mails: khudian@umist.ac.uk, khudian@thsun1.jinr.ru,
khudian@sc2a.unige.ch

   \medskip

  {\bf Abstract:}
   We consider semidensities on a supermanifold $E$ with
   an odd symplectic structure.
     We define a new $\Delta$-operator action on semidensities
     as the proper framework for Batalin-Vilkovisky formalism.
    We establish relations between semidensities on $E$
    and differential forms on Lagrangian surfaces.
    We apply these results to Batalin-Vilkovisky geometry.
    Another application is to $(1.1)$-codimensional surfaces in $E$.
    We construct a kind of pull-back of semidensities to such
     surfaces. This operation and the $\Delta$-operator are used
     for obtaining integral invariants for $(1.1)$-codimensional
     surfaces.

 \bigskip
 \centerline {\it Contents}

\item {1.} Introduction. . . . . . . . . . . . . . . .
         . . . . . . . . . . . . . . . . . . . . .
              . . . .  2
 \item{2.} $\Delta$-operator on semidensities
 . . . . . . . . . . . . . . . . . . . . . . . .
 . . . . . . .  5
 \item {3.} Differential forms on cotangent bundle
 and semidensities
 . . . . . . . . . . . . .  11
 \item {4.}  Semidensities on $E$
   and differential forms on even
    Lagrangian surfaces
    . . . . .  15
  \itemitem {4.1} {\it Identifying symplectomorphisms
   for  even Lagrangian surfaces} . . . . . . 16
    \itemitem {4.2} {\it
   Relation between
  semidensities and differential forms}
   \itemitem {}{\it on a Lagrangian surface}. . . . . . . . . . .
    . . . . . . . . . . . . . . . . . . . 22
\itemitem {4.3} {\it  Application to BV-geometry}
. . . . . . . . . . . . . . . . . . . . . . . . . . . 27
  \item {5.} {Invariant densities on surfaces}
  . . . . . . . . . . . . . . . . . . . . . . . . . . . . . . 29
\item {}  Discussion
. . . . . . . . . . . . . . . . . . . . . . .
 . . . . . . . . . . . . . . . . . . . 35
\item {}  Acknowledgment
. . . . . . . . . . . . . . . . . . . . . . . . . . . . . . . . .
. . . . . .$\,$ 36
  \item {}  Appendix 1.  $\L$-points of supermanifolds
  . . . . . . . . . . . . . . . . . . . .. . . . . $\,$ 37
  \item {} Appendix 2.  A simple proof of the Darboux Theorem
   for odd symplectic structure  38
  \item {} Appendix 3.
  Hamiltonians of adjusted canonical transformations. . . . . . . . .
   . $\,$  41
\item {} References
. . . . . . . . . . . . . . . . . . . . . . . . . . . . . . . . .
. . . . . . . . . $\,${ 43}

\vfill \eject
\medskip
                 \centerline{\bf 1. Introduction}

                  \medskip
\def\superspace{supermanifold $\,$}

   A density of weight $\s$ is a function on a manifold (supermanifold)
   subject to the condition that under change of coordinates it is
   multiplied by the $\s$-th power of the determinant (Berezinian)
    of the transformation. The density of the weight $\s=1$ is
    a volume form.
   (We avoid discussion of orientation here.)

 In this paper
   we study semidensities (densities of weight $\s=1/2$) on
  a  \superspace provided with an odd symplectic structure
   (an odd symplectic supermanifold).
   We introduce a differential operator $\Delta$, which acts on
   semidensities.
  Our considerations
 lead to a straightforward geometrical interpretation
of Batalin-Vilkovisky master equation. On the other hand, we
elaborate a new outlook on the invariant semidensity defined on
$(1.1)$-codimensional surfaces embedded in odd symplectic
supermanifold [\CMP] and construct integral invariants on these
surfaces.

\smallskip

  The concept of an odd symplectic \superspace
  and $\Delta$-operator on it appeared in mathematical physics
 in the pioneer works of I.A.Batalin and G.A.Vilkovisky [\BVa, \BVb],
  where these objects were used for constructing covariant
   Lagrangian version of the BRST quantization (BV formalism).
    The geometrical meaning of these objects
    and interpretation of BV master equation in its terms
    were studied in [\JMPa, \MPL, \JMPd]
  and most notably by A.S.Schwarz in [\SchCMPa].

 Let us shortly sketch the results of [\JMPa, \MPL, \JMPd, \SchCMPa].

     If an odd symplectic
     \superspace is provided with a volume form $d{\bf v}$,
     then one can consider operator
     $\D$ such that its action on a
     function on this \superspace is equal
     (up to a coefficient) to the divergence
      of the Hamiltonian vector field
     corresponding to this  function
     w.r.t. the volume form $d{\bf v}$ [\JMPa].
    This second order differential operator is not trivial because
      transformations
      preserving odd symplectic structure do not preserve
      any volume form (Liouville theorem fails to be fulfilled in
      a case of an odd symplectic structure).

   We call coordinates $z^A=\{x^1,\dots,x^n,\theta_1,\dots,\theta_n\}$
   in an odd symplectic supermanifold Darboux coordinates if
     in these coordinates
     Poisson bracket corresponding to the symplectic structure has
     the canonical form:
     $\{x^i,\theta_j\}=\delta^i_j$,
      $\{x^i,x^j\}=0$.

     Consider a special case, where a volume form in some
     Darboux coordinates is just the coordinate volume form:
                    $$
    \dv=D(x,\theta),\qquad (D(x,\theta)=
    dx^1\dots dx^n d\theta_1\dots d\theta_n)\,.
       \eqno (1.1)
            $$
      In the following we shall refer to it
      as to a particular condition for a volume form.
      Then in this case the operator $\D$
      is given by the following explicit formula
                         $$
                     \D=
                     \sum_{i=1}^n{\p^2 \over \p
                     x^i\p\theta_i}\,,
                           $$
   and it obeys the condition
                             $$
                       \D^2=0\,.
                                           \eqno (1.2)
                                $$
  (See Section 2 for details.)

    The concept of an odd symplectic
    supermanifold provided with a volume form
     is crucial in the geometrical interpretation
      of BV formalism.

    Let $f$ be an even function on
    odd symplectic \superspace
     with a coordinate volume form (1.1) in some Darboux coordinates
    and let $d{\bf v}^\prime=fd{\bf v}$ be
    a new volume form on it.
      In general, for the new volume form  $d{\bf v}^\prime$
       neither condition (1.1) in some Darboux coordinates,
        nor condition (1.2)
       are true.
    The main essence of the geometrical
     formulation of BV formalism
     can be shortly expressed in the following two
        statements [\MPL, \SchCMPa, \JMPd]:

    \smallskip

         {\bf Statement 1.} (see [\MPL, \SchCMPa, \JMPd])

   {\it Consider the following three conditions on
   the volume form $d{\bf v}^\prime=fd{\bf v}$
    and the corresponding $\Delta$-operator:}
                      $$
                      a)\quad
             \matrix
               {
   \hbox{\it there exist Darboux coordinates
       such that the volume form $d{\bf v}^\prime=f\dv$}\cr
   \hbox { {\it has appearence} (1.1) {\it in these
      coordinates}}\,,\cr
                    }
                       \eqno (\statement a)
             $$

             $$
             b)\quad\quad
             \matrix
                 {
          \D\sqrt f=0\,,\cr
    \hbox {\it (BV master-equation for the master-action
     $S=\log\sqrt f$)}\,,\cr
               }
                          \eqno (\statement b)
                      $$

              $$
    c)\quad \Delta^2_{d{\bf v^\prime}}=0\,.
                           \eqno (\statement c)
              $$
{\it Implications}
               $$
a)\Rightarrow b)\Rightarrow c)
               $$
{\it hold. The conditions} a), b), c) {\it are equivalent under
some assumptions
 (see details below)}.

\smallskip

 {\bf Statement 2.} (see [\SchCMPa])

  {\it The integrand of the BV partition function is a semidensity
   $\sqrt{f\dv}$,
   which is a natural integration object over Lagrangian
   surfaces in  odd symplectic supermanifolds.
   In the case if condition} (\statement b) {\it is fulfilled,
    the corresponding integral does not change under
    small variations of Lagrangian surface (gauge-independence
     condition)}.

 The analysis of these statements in [\MPL, \SchCMPa, \JMPd]
  is particularly
based on the following geometrical observations.

 Let $\Pi T^*M$ be the \superspace associated with the cotangent bundle
  $T^*M$ for an arbitrary manifold $M$.
($\Pi T^*M$ is obtained by
   changing the parity of fibers in $T^*M$.)
    Functions on $\Pi T^*M$
   correspond to multivector  fields on $M$.
  The \superspace $\Pi T^*M$ is
    provided with a canonical
   odd symplectic structure.
   The Schouten  bracket of multivector fields
   on it corresponds to the odd Poisson bracket
    of functions on $\PTM$.
  The manifold $M$ is a Lagrangian surface in $\Pi T^*M$.
  If $dv$ is a volume form on $M$
   then the odd sympelctic \superspace $\Pi T^*M$
    provided with a volume form $d{\bf v}=dv^2$
     satisfies conditions (1.1) and (1.2).
    In this case
    the action of operator $\D$ on function on $\PTM$
    corresponds to the
    divergence operator on
 multivector fields on $M$.

   The most profound and detailed analysis of these constructions
    and their relations with Statements 1 and 2 was performed
   in the paper [\SchCMPa]. Particularly in this paper
   some important relations were established between
 differential forms on $M$ and volume forms
 in $\Pi T^*M$ and it was observed
  that the square root of an arbitrary volume
 form in an odd symplectic \superspace is a natural integration object
 over arbitrary Lagrangian surfaces in this supermanifold.

\smallskip

  In this paper we
   consider an odd symplectic supermanifold
   $E=E^{n.n}$. We consider semidensities on $E$.
   We define a new operator $\Delta$ which acts on semidensities.
   Our new operator $\Delta$ is related with operator considered above,
   but it does not require any additional structure on $E$.
  We see that semidensities in an odd symplectic supermanifold,
  not volume forms (densities)
   are naturally related
   with differential forms on even Lagrangian surfaces.
   In particularly the action of $\Delta$-operator on semidensities
     corresponds to the action
     of exterior differential on differential forms.
        The detailed analysis of group of canonical transformations
         for arbitrary odd sympelctic supermanifold $E$
        leads us to establishing relations between calculus of
        semidensities on $E$ and calculus of differential
        forms on even Lagrangian surfaces.

    Our considerations have following two applications.

      In terms of semidensities BV master equation (\statementb)
        gets an invariant formulation and the difference between
         conditions (\statementa, b, c) can be formulated exactly.
         (We note that in papers [\JMPd] and [\SchCMPa] was
      stated that conditions (\statement a), (\statement b) and (\statement c))
      are equivalent, in spite of the fact that difference between these
      conditions was formulated in non-explicitly way in
      Theorem 5 of the paper [\SchCMPa].)

   On the other hand, our considerations give us a new approach
   for obtaining invariant densities and corresponding
   integral invariants on surfaces embedded in
   odd symplectic \superspace with a volume form.
  (The problem of constructing integral invariants for an {\it odd}
  symplectic structure drastically differs from the corresponding problem
  for usual symplectic structure (see in details Section 5)).

  The exposition is organized as follows.

    In Section 2 we recall the basic definitions
   of odd symplectic supermanifold and
   the properties of $\Delta$-operator acting on functions
   provided a volume form is chosen.
 Then we consider semidensities and
 give the intrinsic definition of the
  $\Delta$-operator acting on semidensities
   in an odd symplectic supermanifold. Using this operator
   we formulate  $BV$ master equation
   in invariant way.

  In Section 3 we analyze these objects in terms of
   underlying even geometry considering as the basic
 example
 a supermanifold $\PTM$ associated with the cotangent bundle
 of the usual manifold $M$.
 We establish correspondence between differential forms on $M$
 and semidensities on the supermanifold $\PTM$ and
 analyze the basic formulae of the calculus of differential forms
 in terms
 of semidensities. We also come to new algebraic operations
  on differential
 forms which naturally appears in terms of semidensities.

  In Section 4 we consider even ($(n.0)$-dimensional) Lagrangian
  surfaces
  in odd symplectic supermanifold $E$
  and study correspondence between differential forms on these
  Lagrangian surfaces and semidensities on $E$.
  For any given even Lagrangian surface $L$ this
  correspondence depends on symplectomorphism
  identifying $\Pi T^*L$
  with $E$.
   We prove an existence of an
   identifying symplectomorphism, study identifying symplectomorphisms
   and corresponding subgroups of canonical transformations, and
  investigate in details at what extent
  the correspondence between semidensities and differential forms
  depends on a choice of Lagrangian surface
  and identifying symplectomorphism.
  On the base of these considerations we come to
 statements that generalize
 results of the paper [\SchCMPa] and
 we  formulate exactly differences between conditions
 (\statement a),\statement b) and (\statement c)
  for BV formalism geometry.

  In Section 5
  we provide a natural interpretation
of the odd invariant semidensity on $(1.1)$-codimensional surfaces
that was constructed in [\But, \CMP]. We show that
 this semidensity can be considered as a kind of
 "pull-back" of a semidensity in the
  ambient odd symplectic \superspace.
   This leads us  to construction of another semidensity
   and  two densities (integral invariants), even and odd,
    of rank $k=4$
   on  $(1.1)$-codimensional surfaces.
   These densities seems to be the simplest (having the lowest rank)
   non-trivial integral invariants on surfaces in
  odd symplectic \superspace provided with a volume form.

  The paper contains also three appendices.

 In Appendix 1 we briefly sketch the definition of supermanifold
 as a functor from category of Grassmann algebras to category of sets,
 suggested and elaborated by A.S.Schwarz [\SchCMPb] (see also [\leites]),
 and which we use throughout this paper.
 This definition makes possible to use a language of points
 for supermanifolds. (For basic definitions of constructions in
 supermathematics see books [\Berezin, \leites, \Vor].)

  In Appendix 2 we deliver a simple proof for Darboux theorem
  for odd symplectic supermanifold.

  In Appendix 3 we prove a technical result
  about canonical transformations generated by Hamiltonians.

 \bigskip
 \centerline {\bf 2. $\Delta$-operator on Semidensities}
   \medskip

   In this Section we recall the definitions and properties
    of  odd symplectic \superspace and of the $\Delta$-operator on functions.
   Then we consider semidensities in odd symplectic \superspace
    and define  the action of $\Delta$-operator on semidensities.
    Compared to functions this definition is intrinsic and
    does not require any additional structures (like volume form).

         Let $E^{n.n}$  be $(n.n)$-dimensional \superspace
         and
  \def\Darbouxcoord{$z^A=$  $\{x^1,\dots,x^n$, $\theta_1,\dots,\theta_n\}\,$}
        \Darbouxcoord,
         be local coordinates on it
         ($p(x^i)=0,\, p(\theta_j)=1$, where $p$ is a parity).

        We say that this \superspace is odd symplectic \superspace
       if it is endowed with an odd symplectic structure, i.e.
       an odd  closed non-degenerate $2$-form
           $\Omega=\Omega_{AB}(z)dz^A dz^B$
       ($p(\Omega)=1$, $d\Omega=0$)
  is defined  on it [\Berezin, \leitb, \leites].

  In the same way as in the standard symplectic calculus one can
 relate to the odd symplectic structure the odd Poisson
      bracket (Buttin bracket) [\buttin,\Berezin, \leitb, \leites]:
      \def\bracket {2.1}
                       $$
                     \{f,g\}=
          {\p f\over\p z^A}
                     (-1)^
                   {p(f)p(z^A)+p(z^A)}
                          \Omega^{AB}
              {\p g\over\p z^B}\,,
              \eqno (\bracket)
                       $$
     where $\O^{AB}=\{z^A,z^B\}$
 is the inverse matrix to $\O_{AB}$ :
$\O^{AC}{ \O_{CB}}=\d^A_B$ .

  Hamiltonian vector field
  \def\hamvectorfield {2.2}
                       $$
         {\bf D}_f=\{f,z^A\}
      {\p\over\p z^A}\,,
                 \quad
         {\bf D}_f(g)=\{f,g\},\quad
          \O({\bf D}_f,{\bf D}_g)=-\{f,g\}
                                     \eqno (\hamvectorfield)
                      $$
corresponds to every function $f$.

    The condition of the closedness of the form defining
     symplectic structure implies the Jacoby
    identity:
          \def\jacoby {2.3}
                  $$
      \def\jacoby {2.3}
                  \{f,\{g,h\} \}(-1)^
              {(p(f)+1)(p(h)+1)}+\hbox{cycl. permutations} =0\,.
                      \eqno (\jacoby)
              $$.

    Using the analog of the Darboux Theorem [\Shander, \SchCMPa]
   (see also Appendix 2) one can consider in
   a vicinity of an arbitrary point
   coordinates \Darbouxcoord  such that in
   these coordinates  symplectic structure  and the
   corresponding odd Poisson bracket
    have locally the canonical expressions
    \def\darbouxtheorem {2.4}
                          $$
            \Omega=I_{AB}dz^A dz^B \colon\quad
            \Omega\left(
     {\partial\over\partial x^i},
     {\partial\over\partial x^j}
            \right)=0,\,
           \Omega\left(
     {\partial\over\partial \theta^i},
 {\partial\over\partial \theta^j}\right)=0,\,
 \Omega\left({\partial\over\partial x^i},
 {\partial\over\partial \theta^j}\right)=-\delta_{ij}\,,
                      $$
    and respectively
                        $$
         \{x^i,x^j\}=0,\, \{\theta_i,\theta_j\}= 0,\,
                 \{x^i,\theta_j\}=-\{\theta_j,x^i\}=\d^i_j\,,
                          $$
                          $$
                      \{f,g\}=
                      \sum_{i=1}^n
                       \left(
         {\p f\over\p x^i}
         {\p g\over\p \theta_i}
                           +(-1)^{p(f)}
         {\p f\over\p \theta_i}
         {\p g\over\p x^i}
                       \right)\,.
                                           \eqno (\darbouxtheorem)
                         $$
   These coordinates are called Darboux coordinates.
  Transformation of Darboux coordinates to another Darboux coordinates
 is called canonical transformation of coordinates. Respectively
transformation of supermanifold that transforms Darboux
coordinates
 to another Darboux coordinates is called canonical
 transformation.

   We consider also odd symplectic supermanifold provided
   additionally with a volume form:
                         $$
     d{\bf v}=\rho(z) Dz=\rho(x,\theta)D(x,\theta)\,,\quad
                  (p(\rho)=0)\,.
                                     \eqno (2.5)
            $$
  $Dz=D(x,\theta)$ is coordinate volume form
  ($D(x,\theta)=dx^1...dx^n d\theta_1\dots d\theta_n$).
  Coordinate volume forms in different coordinates
  are related by Berezinian (superdeterminant) of
  coordinate transformation [\Berezin]:
  \def\berdef{2.6}
                     $$
         {D\t z\over Dz}={\rm Ber}{\p\t z\over\p z},
         \quad {\rm where}\,\,
                      {\rm  Ber}
                    \pmatrix
                 {I_{00}&I_{01}\cr I_{10}&I_{11}\cr}
                        =
           {\det\, (I_{00}-I_{01}I_{11}^{-1}I_{10})\over \det\, I_{11}}\,.
                     \eqno (\berdef) 
                     $$

We suppose that the volume form (2.5) is non-degenerate, i.e. for
the every point $z_0$ the number part $m(\rho(z_0))$ of
$\rho(z_0)$ is not equal to zero.

\smallskip

In the paper [\JMPa] we show that an odd symplectic structure (in
fact an odd Poisson bracket structure, which might be degenerate)
and a volume form allow to define the $\Delta$-operator (or
Batalin-Vilkovisky operator; this is the invariant formulation of
the operator introduced in BV-formalism [\BVa]).
 The construction is as follows.
The action of $\Delta$-operator on an arbitrary function in an odd
symplectic supermanifold provided with a volume form is equal (up
to coefficient) to the divergence w.r.t. volume form (2.5) of the
Hamiltonian vector field corresponding to this function.
 Using (\hamvectorfield) we come to  the formula
\def\deltadefgen{2.7}
                              $$
    \D f={1\over 2}(-1)^{p(f)}{\rm div}_{d{\bf v}}{\bf D}_{f}=
             {1\over 2}(-1)^{f}
         \left(
         (-1)^{p({\bf D}_{f}z^A+z^A)}
           {\p\over\p z^A}
             \{f,z^A\}+
               D_f^A{\p \log\rho(z)\over\p z^A}
                    \right)\,.
                       \eqno (\deltadefgen) 
                          $$

 In Darboux coordinates:
 \def\deltadef {2.8}
                        $$
    \D f=\Delta_0 f+ {1\over 2}\{\log\rho, f\}\,,
                   \eqno (\deltadef) 
            $$
 where $\rho(z)$ is given by (2.5) and
 \def\deltatrdef{2.9}
                   $$
     \Delta_0 f=\sum_{i=1}^n{\p^2 f\over \p x^i\p\theta_i}\,.
                         \eqno (\deltatrdef) 
                   $$

$\Delta$-operator on functions satisfies the relations [\BVb,
\MPL] :
\def\deltaprop {2.10}
            $$
    \D \{f,g\}=\{\D f,g\}+(-1)^{p(f)+1}\{f,\D g\}\,,
                        $$
            $$
    \D (f\cdot g)=\D f\cdot g+(-1)^{p(f)} f\cdot\D g+
    (-1)^{p(f)}\{f,g\}\,.
                                \eqno (\deltaprop)
                        $$

Operator $\Delta_0$ in (\deltatrdef) is not an invariant operator
on functions (i.e. it depends on the choice of Darboux
coordinates). It can be considered as $\D$ operator for
coordinate volume form $D(x,\theta)$ in the chosen Darboux
coordinates
 \Darbouxcoord. If
 $\t z^A=$ $\{\t
x^1,\dots,\t x^n,$ $\t\theta_1,\dots,\t\theta_n\}$ are another
Darboux coordinates then from (\deltadef) it follows that
\def\deltatransform {2.11}
                           $$
                           \Delta_0 f=
                           \widetilde\Delta_0 f+
                           {1\over 2}
                           \{\log {\rm Ber}{\p z\over\p\t z},f\}\,,
   \eqno (\deltatransform) 
                            $$
where $\widetilde\Delta_0$ is operator (\deltatrdef) in Darboux
coordinates $\t z^A=\{\t x^1,\dots,\t
x^n,\t\theta_1,\dots,\t\theta_n\}$.

\smallskip

   Now we consider semidensities
   on odd symplectic supermanifold.
   In local coordinates $z^A=\{x^i,\theta_j\}$ they have the appearance
   ${\bf s}=s(z)\sqrt{Dz}=$ $s(x,\theta)\sqrt{D(x,\theta)}$.
 Under coordinate transformation
 $z^A=z^A(\t z)$ the coefficient $s(z)$
 is multiplied by the square root
 of the Berezinian of corresponding transformation:
$s(z)\mapsto s(z(\t z)){\rm Ber}^{1/2}(\p z/\p\t z)$.

 We shall define a new operator, which we denote $\Ds$,
 and which will act on the space of semidensities.

\smallskip

 {\bf Definition}
{\it Let ${\bf s}$ be a semidensity
 and $s(z)\sqrt{Dz}$ be its local expression
in some Darboux coordinates $z^A=\{x^1,\dots, x^n,$
$\theta_1,\dots,\theta_n\}$. The local
 expression for the semidensity $\Ds {\bf s}$
 in these coordinates is given by the following formula:}
\def\deltadensdef {2.12}
                     $$
              \Ds {\bf s}=(\Delta_0 s(z))\sqrt{Dz}=
     \sum_{i=1}^n{\p^2 s\over \p x^i\p\theta_i}
              \sqrt{D(x,\theta)}\,.
                                  \eqno (\deltadensdef) 
                 $$
 {\it The semidensity $\Ds{\bf s}$ is an odd (even) if semidensity
  ${\bf s}$ is an even (odd) semidensity, thus $\Ds$ is an odd operator.}

Contrary to the operator $\D$ on functions, the operator $\Ds$ on
semidensities does not need any volume structure.

  \smallskip

To prove that $\Ds$-operator is well-defined by  formula
(\deltadensdef),
 one has to check that r.h.s. of
 (\deltadensdef) indeed defines
semidensity, i.e. if $\{\t z^A\}=\{\t x^1,\dots,\t
x^n,\t\theta_1,\dots,\t\theta_n\}$ are another
 Darboux coordinates then
\def\correctness {2.13}
                   $$
           \left(
  \sum_{i=1}^n{\p^2 \over \p x^i\p\theta_i}
                 s(z)
         \right)_{z(\t z)}
                  \cdot
            \left(
     {\rm Ber}{\p z(\t z)\over\p\t z}
            \right)^{1/2}
                    =
               \sum_{i=1}^n
      {\p^2\over \p \t x^i\p\t\theta_i}
              \left(
                  s(z(\t z))
                    \cdot
      {\rm Ber}^{1/2}{\p z(\t z)\over\p \t z}
                      \right)\,.
                \eqno (\correctness) 
             $$
    First of all we check this condition for infinitesimal
     canonical transformations.
   They
  are generated by an odd
  function (Hamiltonian) via the corresponding Hamiltonian vector filed.
  To an odd Hamiltonian $Q(z)$ corresponds infinitesimal
  canonical transformation
   $\t z^a=z^A+\varepsilon \{Q,z^A\}$ generated by
   the vector field ${\bf D}_Q$ in (\hamvectorfield).
 To the action of this transformation on the semidensity ${\bf s}$
     corresponds differential
    $\d_Q(s\sqrt{Dz})=$ $\Delta_0 Q\cdot s\sqrt{Dz}-$
    $\{Q,s\}\sqrt{Dz}$,
    because $\d s=-\varepsilon\{Q,s\}$ and
    $\d Dz=\varepsilon\d {\rm Ber}(\p z/\p\t z)Dz=
    2\Delta_0Q Dz$.
    Using that $\Delta_0^2=0$ and relation (\deltaprop)
     we come to commutation relations
   $\Delta_0\d_Q=\d_Q\Delta_0$.
     Thus we come to condition (\correctness),
    for infinitesimal transformations.

  To check condition  (\correctness)
  for arbitrary canonical transformation
  we need the following

  {\bf Lemma 1}
\def\cantransform {2.14}
 1. {\it Every canonical transformation of Darboux coordinates
 $\t z^A=\F^A(z)$
  can be decomposed into canonical transformations
  $\F(z)=\F_{s}\left(\F_p\left(\F_{\rm adj}(z)\right)\right)$,
   where

  a) canonical transformation $\t z=\F_{\rm adj}(z)$,
     has the following form}
                                    $$
                                    \cases
                                    {
     \t x^i(x,\theta)\big\vert_{\theta=0}=x^i\,,\cr
     \t\theta_i(x,\theta)\big\vert_{\theta=0}=0\,,\cr
                          }\qquad
       (i=1,\dots,n)\,,
                                          \eqno (\cantransform a) 
                                     $$
  {\it we call later this canonical transformation
  of Darboux coordinates adjusted
  canonical transformation;

   b) canonical transformation
  $\t z=\F_p(z)$ has the form}
                                      $$
                                   \cases
                                     {
                                  \t x^i=x^i(x)\cr
                                  \t \theta_i=
                         {\p  x^m(\t x)\over\p \t x i}\theta_m\cr
                                     }\,\qquad
                                (i,j=1,\dots,n)\,,
                                     \eqno (\cantransform b) 
                                    $$
  {\it we call later this canonical transformation
  of Darboux coordinates, which is generated
    by transformation $\t x^i=x^i(x)$
    "point"-canonical transformation;

   c) canonical transformation
   $\t z=\F_s(z)$ has the following form}
                                     $$
                                     \cases
                                       {
           \t x^i=x^i\cr
           \t\theta_i=\theta_i+\Psi_i( x)\,\cr
                        }\quad
                        {\rm such\,that}\,\,
   {\p\Psi_i( x)\over \p \ x^j}-{\p\Psi_j(x)\over \p x^i}=0\,,
                         (i=1,\dots,n)\,,
                                          \eqno (\cantransform c) 
                                  $$

  {\it we call later this canonical transformation
  of Darboux coordinates special canonical transformation.}

  2.  {\it Berezinian
  of adjusted canonical transformation} (\cantransform a)
   {\it obeys to the condition
    ${\rm Ber}{\p \t z\over \p z}\big\vert_{\theta=0}=1$,
    Berezinian of "point"
  canonical transformation} (\cantransform b)
   {\it is equal to $\det^2{\p \t x\over\p x}$,
    and Berezinian of special canonical transformation}
        (\cantransform c)
    {\it is equal to one.

    In particular, numerical part of Berezinian
  of arbitrary canonical transformation is positive.}

  3. {\it Adjusted canonical and special canonical transformations
  of Darboux coordinates}
  (\cantransform a, \cantransform c)
  {\it are canonical transformations generated by Hamiltonian, i.e.
      they  can be included in
      one-parametric family of canonical transformations
     of Darboux coordinates generated by an odd Hamiltonian}:
\def\hamiltoniandeformation {2.15}
     $$
  \exists Q(z,t)\colon\quad
     {dz_t\over dt}=\{Q,z_t\},\,\, 0\leq t\leq 1\,\, {\rm such\, that}\,\,
      z_0=z\,, z_1=\t z\,.
      \eqno (\hamiltoniandeformation) 
      $$
{\it Special canonical transformation} (\cantransform c) {\it is
generated locally} {\it by Hamiltonian $Q=Q(x)$, such that $\p_i
Q(x)=\Psi_i(x)$}.
 {\it There exists unique "time"-independent Hamiltonian
    $Q=Q(z)$ obeying to condition
    $Q(x,\theta)=Q^{ik}\theta_i\theta_k+\dots$, i.e.
  $Q=O(\theta^2)$
 that generates}
 {\it given adjusted canonical transformation.}

   \smallskip

    Prove this Lemma.

   Let $\t z=\F(z)$ be arbitrary canonical transformation:
   $\t x^i=f^i(x,\theta)=f_0^i(x)+O(\theta)$ and
   $\t \theta_i=\Psi_i(x)+O(\theta)$.
    Consider coordinates $\{{\bar z}^A\}=$
    $\{{\bar x^i}, {\bar\theta_i}\}$
    that are related with coordinates $\{\t z^A\}$
    by the following special canonical transformation:
    $\t z^A=\F_s({\bar z})$ such that
    $\t x^i={\bar x^i}$ and
    $\t \theta_i={\bar\theta_i}+\Psi_i(g({\bar x}))$,
    where $g\circ f_0={\bf id}$. Then
    ${\bar x}^i=f^i(x,\theta)$ and
   ${\bar \theta}_i=O(\theta)$. Now
   consider coordinates $\{z^{\prime A}\}=$
   $\{x^{\prime i}, \theta_i^\prime\}$
    that are related with coordinates $\{\bar z^A\}$
    by the  "point" canonical transformation
    $\bar z^A=\F_p({z^\prime})$
    generated by functions $\bar x^i=f_0^i(x^\prime)$.
    Then it is easy to see that
    initial coordinates  $\{z^A\}$
    are related with coordinates
    $\{z^{\prime A}\}$ by adjusted canonical transformation
    $z^{\prime A}=\F_{\rm adj}(z)$:
    $x^{\prime i}=x^i+O(\theta)$, $\theta^\prime_i=O(\theta)$.

  The second statement of Lemma can be proved by easy straightforward
  calculation of Berezinian (\berdef) for transformations
  (\cantransform a, \cantransform b, \cantransform c).

    We perform the proof of the statement 3 of
    Lemma for adjusted canonical
    transformations in Appendix 3. \finish

   \smallskip

   Now we return to the proof of relation (\correctness).

   First we note
   that from second statement of Lemma it follows that
    square root operation in (\correctness) is well-defined.

   From decomposition (\cantransform) it follows that
    it is sufficient to check condition (\correctness)
   separately for adjusted, "point", and special
   canonical transformations.
   From the third statement of Lemma it follows that
   for adjusted and special  canonical
   transformations  the condition
   (\correctness) can be checked only
   infinitesimally and this is performed already.
   For "point" canonical transformation (\cantransform b)
   the condition (\correctness)
    can be easily checked
    strai\-ght\-forwardly using (\deltaprop), (\deltatransform)
    and the fact that Berezinian of this transformation
    does not depend on $\theta$. \finish

\smallskip

    The action of differential $\d_Q$ corresponding to
    infinitesimal canonical transformation
    on semidensities can be rewritten in an explicitly invariant way:
\def \denstransform {2.16}
                  $$
          \d_Q{\bf s}=Q\cdot\Ds{\bf s}+\Ds( Q{\bf s})=
                [Q,\Ds]_+{\bf s}\,.
                          \eqno (\denstransform) 
                        $$

On an odd symplectic supermanifold provided with a volume form
$\dv$ (density of the weight $\sigma=1$)
 we can construct new invariant objects, expressing them via the
semidensity related with volume form and operator $\Ds$:
\def\denscalculus {2.17}
                       $$
        {\bf s}=\sqrt\dv\qquad
      \hbox {semidensity ($\sigma={1\over 2}$)}\,,
                                   \eqno (\denscalculus a) 
            $$
            $$
        \Ds {\bf s}=\Ds\sqrt\dv\qquad
        \hbox {semidensity ($\sigma={1\over 2}$)}\,,
                        \eqno (\denscalculus b) 
            $$
             $$
    {\bf s}\Ds{\bf s}=\sqrt\dv\Ds\sqrt\dv\qquad
    \hbox {density ($\sigma=1$)}\,,
                           \eqno (\denscalculus c)
            $$
            $$
             {1\over{\bf s}}\Ds{\bf s}=
         {1\over\sqrt\dv}\Ds\sqrt\dv\qquad
         \hbox {function ($\sigma=0$)}\,.
                                  \eqno (\denscalculus d)
             $$

  From definition (\deltadensdef) of $\Ds$-operator and relations
  (\deltaprop) it follows that
$\Ds$-operator obeys to the following properties:
\def\deltadenspropa {2.18}
                      $$
              (\Ds)^2=0\,,
               $$
                      $$
            \Ds (f\cdot \sqrt\dv)=
            (\D f)\cdot \sqrt\dv+
            (-1)^{f}f\cdot \Ds \sqrt\dv\,,
                                \eqno (\deltadenspropa) 
                    $$
  and
\def\deltadenspropb {2.19}
                 $$
             \D^2 f=
         \{{1\over\sqrt\dv}\Ds\sqrt\dv,f\}\,.
                           \eqno (\deltadenspropb) 
            $$

  We call semidensity
  $\bf s$ {\it closed} semidensity if $\Ds {\bf s}=0$
 and we call ${\bf s}$ an {\it exact} if  there exists another
 semidensity ${\bf r}$ such that ${\bf s}=\Ds {\bf r}$.

   In the case if an odd symplectic supermanifold
    is provided with a volume form $\dv$
    such that this volume form is equal to coordinate volume form
    $D(\t x,\t\theta)$
    in some Darboux coordinates $\{\t x^1,\dots,\t x^n$,
    $\t\theta_1,\dots,\t\theta_n\}$
   then evidently
\def\batsemidensity {2.20}
$\Ds \sqrt {\dv}=0$. Considering this relation in
 another Darboux coordinates \Darbouxcoord
we come to the formula
                     $$
     \Delta_0{\rm Ber}^{1/2}\left({\p \t z\over \p z}\right)=
         \sum_{i=1}^n{\p^2 \over \p x^i\p\theta_i}
               {\rm Ber}^{1/2}
               \left(
               {\p (\t x,\t\theta)\over \p (x,\theta)}
               \right)=0\,.
                        \eqno (\batsemidensity) 
                     $$
We note that formulae (\deltatrdef) and (\deltatransform)
 for $\Delta_0$ operator
 were first studied
 by I.A.Batalin and G.A.Vilkovisky ([\BVa, \BVb]).
 In particular they obtained formula (\batsemidensity).
These results receive its clear geometrical interpretation
 in terms of semidensities and action of $\Ds$-operator on them.

 We say that semidensity ${\bf s}=s(x,\theta)\sqrt{D(x,\theta)}$ is
 non-degenerate if a number part $m(s(x,\theta))$
 of $s(x,\theta)$ is not equal
 to zero at any $x$.
  Every volume form defines non-degenerate even semidensity
  by relation (\denscalculus a)
  and respectively volume form corresponds to
  every non-degenerate even semidensity.
   We say that even non-degenerate semidensity ${\bf s}$ obeys to
   BV-master equation if it is closed
   and we denote by $\B_{\rm deg}$ a set of these densities.
\def\bvsdensdef{2.21}
                     $$
        \B_{\rm deg}=
        \{{\bf s}\colon\quad
        \Ds{\bf s}=0\,,\quad p(s(x,\theta))=0,\quad
            m(s(x,\theta))\not=0\}\,.
        \eqno (\bvsdensdef) 
                 $$
BV-master equation (condition (\statementb))
 was not formulated invariantly in
[\MPL,\SchCMPa]. Condition $\Ds {\bf s}=0$ (closedness of
semidensity ${\bf s}$) gives invariant formulation to BV
master-equation.

\bigskip

      \centerline {\bf 3. Differential forms on
      cotangent bundle and semidensities}

\medskip

 We consider in this section basic example of an odd symplectic supermanifold
 yielded by cotangent bundle of usual manifolds.
 We clarify geometrical meaning of previous constructions
 and establish relations
 between differential forms on manifold and semidensities
 on this odd symplectic supermanifold.

  \smallskip

  In the standard sympelctic calculus cotangent bundle of any
  manifold can be provided
  with canonical symplectic
   structure and  it can be considered as basic example of
   symplectic manifold [\Gui].
  Basic example of an odd symplectic supermanifold
  is constructed in the following way.
  Let $M$ be arbitrary $n$-dimensional manifold and
 $T^*M$ be its cotangent bundle.
  Consider a supermanifold
   $\Pi T^*M$ associated with cotangent bundle
 $T^*M$, changing the parity of fibers of cotangent bundle $T^*M$.
  Let $\{x^1,\dots,x^n,p_1,\dots,p_n\}$ be canonical coordinates
 on $T^*M$ corresponding to arbitrary local coordinates
  $\{x^1,\dots,x^n\}$ on $M$, i.e. for a form $w\in T^*M$
  $p_i(w)=w({\p\over \p x^i})$.
  Canonical coordinates
  $z^A=(x^1,\dots,x^n,\theta_1,\dots,\theta_n)$
  on $\PTM$, ($p(\theta_i)=1$) correspond
  to the canonical coordinates $\{x^1,\dots,x^n,p_1,\dots,p_n\}$
  on $T^*M$.
  Odd coordinates $\{\theta_1,\dots,\theta_n\}$ transform
  via the differential of corresponding
   transformation of coordinates $\{x^i\}$ of underlying space $M$
   in the same way as coordinates $\{p_1,\dots,p_n\}$ on $T^*M$:
                $$
\t x^i=\t x^i(x),\quad\t\theta_i=\sum_{k=1}^n
      {\p x^k(\t x)\over\p \t x^i}\theta_k\,,\quad
           (i=1,\dots,n)\,.
                     \eqno (3.1)
            $$
 We define canonical odd symplectic structure on $\PTM$
   considering these coordinates as Darboux coordinates (\darbouxtheorem).
  Thus we assign
   to every atlas $\left[\{x^i_{(\a)}\}\right]$ of coordinates
   on manifold $M$ an atlas
   $\left[\{x^i_{(\a)},\theta_{j(\a)}\}\right]$
   of Darboux coordinates on supermanifold $\PTM$.
   Pasting formulae (3.1) ensure us that this
   canonical symplectic structure is well-defined.

  Later on we call Darboux
  coordinates
  (3.1) on $\Pi T^*M$ induced by coordinates
  on $M$
    {\it Darboux coordinates adjusted to cotangent bundle structure}.
  Unless otherwise stated we assume further
   that Darboux coordinates in a supermanifold associated with
   cotangent bundle are Darboux coordinates adjusted
   to cotangent bundle structure.
   (Canonical transformations (3.1), induced by coordinate transformations
   on the manifold $M$ are "point" canonical
   transformations (\cantransform b).)

The relations between the cotangent bundle structure on $T^*M$
and the odd canonical symplectic structure on $\Pi T^*M$ reveal
in the properties of the following canonical map $\tau_{_M}$
 between multivector fields on $M$ and
 functions on $\Pi T^*M$:
               $$
           \tau_{_M}
           \left(
           T^{i_1\dots i_k}
           {\p\over\p x^{i_1}}
           \wedge\dots\wedge
           {\p\over\p x^{i_k}}
               \right)=
           T^{i_1\dots i_k}
           \theta_{i_1}\dots\theta_{i_k}\,.
                       \eqno (3.2)
                  $$
 This map transforms the
 Schoutten bracket of multivector fields to the
 odd canonical Poisson
 bracket (Buttin bracket) (2.2)
  of corresponding functions [\leites,\buttin]:
                        $$
            \tau_{_M}
            \left(
            \left[
            {\bf T_1},{\bf T_2}
            \right]
              \right)=
            \{
        \tau_{_M} \left({\bf T_1}\right),
        \tau_{_M} \left({\bf T_2}\right)
                  \}\,.
                       \eqno (3.3)
             $$

    Now we construct a map, that establishes correspondence between
    differential forms on $M$ and semidensities on $\Pi T^*M$.
   We consider arbitrary Darboux coordinates
   $\{x^1,\dots,x^n,$ $\theta_1,\dots,\theta_n\}$
    on $\Pi T^*M$
    adjusted to cotangent bundle structure
    and define this map in these Darboux coordinates
    in the following way:
\def\diffdenscorr {3.4}
                  $$
     \eqalign
         {
    \tau_{_M}^{\#}(1)&=\theta_1\dots\theta_n\sqrt{D(x,\theta)},\cr
    \tau_{_M}^{\#}(dx^i)&=(-1)^{i+1}
 \theta_1\dots
 \widehat\theta_i\dots\theta_n\sqrt{D(x,\theta)}\,,\cr
\tau_{_M}^{\#}(dx^i\wedge dx^j)&=(-1)^{i+j }
    \theta_1\dots
    \widehat\theta_i\dots\widehat\theta_j
    \dots\theta_n\sqrt{D(x,\theta)}\,,\,(i<j),\cr
        &\dots\cr
  \tau_{_M}^{\#}(dx^{i_1}\wedge\dots\wedge dx^{i_k})&=
           (-1)^{i_1+\dots+i_k+k}
      \theta_1\dots
    \widehat\theta_{i_1}\dots\widehat\theta_{i_k}
    \dots\theta_n\sqrt{D(x,\theta)},\,
     (i_1<\dots <i_k), \cr
                 }
        $$
                 $$
     \tau_{_M}^{\#}(f(x)w)=f(x)\tau_{_M}^{\#}(w)\,,\quad
  \hbox{for every function $f(x)$ on $M$}\,,
                            \eqno (\diffdenscorr) 
                               $$
where the sign $\,\,\widehat {}\,\,$
 means the omitting of corresponding term.
   For example
if $M$ is two-dimensional space, then
$\tau^{\#}_{_M}(f(x))=f(x)\theta_1\theta_2\sqrt{D(x,\theta)}$,
$\tau^{\#}_{_M}(w_1dx^1+w_2(x)dx^2)=$

\noindent $(w_1\theta_2-w_2(x)\theta_1)\sqrt{D(x,\theta)}$,
$\tau^{\#}_{_M}(wdx^1\wedge dx^2)=-w\sqrt{D(x,\theta)}$.

  One can rewrite (3.4) in a more compressed way:
                       $$
   \tau^\#_M(w)=\left(
   \int w(x,\xi)\exp(\theta_i\xi^i)d^n\xi
           \right)
           \sqrt{D(x,\theta)}\,,
                                  \eqno (\diffdenscorr a) 
                        $$
where $w(x,\xi)$ is a function corresponding to differential form
$w$ in the supermanifold $\Pi TM$ associated
 to the tangent bundle $TM$:
 $w(x,\xi)=w_{i_1\dots i_k}\xi^{i_1}\dots \xi^{i_k}$.
 Odd coordinates $\{\xi^i\}$
 of the fibers in $\Pi TM$ transform as differentials $\{dx^i\}$:
 $\t x^i=\t x^i(x)$ $\mapsto$
 $\t\xi^i={\p\t x^i\over\p x^k}\xi^k$.
 The square of the map (\diffdenscorr a)
$w\rightarrow (\tau^{\#}(w))^2$ transforms differential forms on
$M$ to density (volume form) on $\Pi T^*M$ and  this map was
considered in [\SchCMPa].

To prove that  (\diffdenscorr) is well-defined  for an arbitrary
Darboux coordinates adjusted to cotangent bundle structure we
note that under arbitrary coordinate transformation (3.1) the
integral in r.h.s. of (3.4a)
 is multiplied on
the $\det (\p \t x/\p x)$ and coordinate volume form
$\sqrt{D(x,\theta)}$
 is divided on the module of this determinant, because for transformation
 (3.1):
\def\bercanonical{3.5}
                   $$
                  {\rm Ber}^{1/2}
                          \left(
             {\p(\t x,\t\theta)\over\p( x,\theta)}
                         \right)
                            =
                   {\rm Ber}^{1/2}
                      \pmatrix
                         {
                   &{\p \t x^i(\t x)\over \p x^k}
                   &{\p \t x^r\over \p x^k}
                     {\p^2  x^m\over \p \t x^r\p \t x^i}\theta_m \cr
                    & 0
                    &{\p  x^k\over \p\t x^i}\cr
                               }
                             =
                             \big\vert\det
                        \left(
                        {\p \t x^i( x)\over \p x^k}
                            \right)\big\vert
                            \,.
                                            \eqno (\bercanonical) 
                            $$

{\bf Remark} Map (\diffdenscorr) establishes
 correspondence only up to a sign factor
because r.h.s. of (3.5) is positive for canonical transformation
induced by any coordinate transformation of $M$.
 In the case if $M$ is orientable manifold
considering only Darboux coordinates such that Jacobian of
coordinate transformations is positive one comes to globally
defined map. Sign factor depends on a orientaion of M.

\smallskip

We say that semidensity ${\bf s}$ corresponds to
 differential form $w$ (to the linear combination of
  differential forms $\sum w_k$) if
  ${\bf s}=\tau^{\#}_M (w)=\tau^{\#}_M (\sum w_k)$.

  In the end of this section
   using correspondence between semidensities
   and differential forms
  we consider some standard constructions of differential
  forms calculus in terms of semidensities and
  two geometrical operations on differential forms, which are
  naturally arisen in terms of semidensities.

   \smallskip

  1) From (\diffdenscorr a) it is easy to see that
  an action of operator $\xi^i{\p\over \p x^i}$ on the function
 $w(x,\xi)$ corresponds to the action
  of exterior differential $d$ on differential form
  and to
  the action of $\Ds$-operator on semidensity, i.e.
  the action of $\Ds$-operator
 corresponds to the action of exterior differential:
 \def\innerdelta {3.6}
                            $$
    \Ds\circ\tau_{_M}^{\#}=\tau_{_M}^{\#}\circ d\,.
                         \eqno (\innerdelta) 
            $$
 Closed (exact) semidensity corresponds
 to closed (exact) differential form.

\smallskip

 2) If the semidensity ${\bf s}$ in $\Pi T^*M$ corresponds to
 volume form (differential top-degree form $w$ on $M$) and
 an odd sympelctic supermanifold $\Pi T^*M$ is provided with volume
 form  such that it is equal to the
 square of this semidensity
  then  the action of operator $\D$ corresponds
  to the divergence w.r.t. to the
 volume form $w$ on $M$:
\def\divdelta {3.7}
              $$
              \D\circ\tau_{_M}=
 \tau_{_M}\circ {\rm div}_w\,\quad {\rm if}\quad \dv={\bf s}^2\,{\rm and}\,\,
                         {\bf s}=\tau^\#_M w\,.
              \eqno (\divdelta)
      $$
 (See also [\JMPd, \SchCMPa].)

\smallskip

3) From (3.2) and (3.4) it follows that
\def\innerdeltab {3.8}
                       $$
             \tau_{_M}^{\#}({\bf T}{\cal c} w)=
    \tau_{_M}({\bf T})\cdot\tau_{_M}^{\#}(w)\,,
             \eqno (\innerdeltab) 
                         $$
where ${\bf T}{\cal c} w$ is the inner product of multivector
field ${\bf T}$ with differential form $w$.

\smallskip

4) The meaning of relation (\denstransform) in terms of
differential forms is following.
 In the special case if Hamiltonian
$Q$ corresponds
 to vector field $T^i(x){\p\over\p x^i}$
 ($Q=T^i(x)\theta_i$), then this Hamiltonian
 induces infinitesimal canonical
 transformation that corresponds to the
 infinitesimal transformation of $M$ induced by the vector field
 $T^i(x){\p\over\p x^i}$.
 From (\innerdelta, \innerdeltab) it follows that in this case the
  standard formula for Lie derivative
 of differential forms
 (${\cal L}_T w=$ $dw{\cal c}T+d(w{\cal c}T)$)
 corresponds to relation (\denstransform).
 In a general case canonical transformations
  of $\Pi T^*M$ destroy cotangent bundle structure and
  mix forms of different degrees.
 For example if
  we consider the action of Hamiltonian
  $Q=L\theta_1\dots\theta_n$
  on a semidensity corresponding to form $w=dx^1\wedge\dots\wedge dx^n$
   then we obtain using (\denstransform) that
  $\d w=dL$.

   \smallskip

   5)  If $a=a_i(x)dx^i$  is $1$-form on $M$ then one can see  that
   \def\dopodin{3.9}
                 $$
         \tau^\#(a\wedge w)=
         a_i{\p s\over\p\theta_i}\sqrt{D(x,\theta)},\quad
                {\rm where}\quad
                \tau^\#(w)={\bf s}\,.
                   \eqno (\dopodin)
                $$
   It is more natural from point of view of
   semidensities to consider the following relation
   between $1$-forms on $M$ and semidensities in $\PTM$.
  Let $a=a_idx^i$ be an odd-valued one-form on $M$ with coefficients in
  arbitrary Grassmann algebra $\L$ (see Appendix 1).
   For this form and arbitrary semidensity
   ${\bf s}=s(x,\theta)\sqrt{D(x,\theta)}$
   consider a new semidensity ${\bf s}^\prime$, which we denote by
   $a\,{\cal d}\,{\bf s}$ such that it is given by relation
   \def\firstoperation {3.10}
                            $$
          {\bf s}^\prime=
          a\,{\cal d}\,{\bf s}=s(x,\theta_i+a_i)\sqrt{D(x,\theta)}\,.
                       \eqno (\firstoperation) 
                   $$
  Respectively if semidensity ${\bf s}$
  corresponds to  differential form $w=\sum w_k$ then we denote by
    $a\,{\cal d}\,w$ differential form such that
  semidensity $a\,{\cal d}\,{\bf s}$
    corresponds to $a\,{\cal d}\,w$.
   From (\firstoperation) and (\diffdenscorr a) it follows that
 \def\firstoperationa {3.11}
                           $$
            a\,{\cal d}\,w=\sum_{p=0}^k {1\over p!}
    \underbrace {a\wedge\dots\wedge a}_{p\quad{\rm times}}\wedge
                          w_{k-p}\,,\quad (k=0,\dots,n)\,.
                          \eqno (\firstoperationa) 
                     $$

 Relations (\firstoperation) and (\firstoperationa)
  define an action of abelian supergroup
 of differential odd valued one-forms on semidensities and
 differential forms.
\smallskip

 6) Consider also the following algebraic operation on differential
  forms that seems very natural from the point of view of
  semidensity calculus. Let $w=\sum w_k$ and $w^\prime=\sum w^\prime_k$
  be differential forms on $M^n$ such that top-degree forms $w_n$
  and $w^\prime_n$ are not equal to zero.
  Then we consider a new form
  \def\secondoperation {3.12}
                            $$
              \tilde w=w*w^\prime\colon\quad
                  \tau^\# \tilde w=
                  \sqrt{\tau_M^\#(w_1)\cdot\tau_M^\#(w_2)}\,.
                     \eqno (\secondoperation) 
                      $$
The condition $w_n\not=0, w^\prime_n\not=0$ for top-degree forms
makes well-defined a square root operation on corresponding
semidensities.

\bigskip

       \centerline {\bf 4.Semidensities on $E$ and differential
       forms  on even Lagrangian surfaces}

   \medskip

In the previous Section we analyzed relations between differential
forms on manifold $M$ and semidensities on supermanifold $\PTM$
using Darboux coordinates in $\PTM$ that are adjusted to
cotangent bundle structure of $T^*M$. (Relations (\diffdenscorr)
are not invariant with respect to an arbitrary canonical
transformation of Darboux coordinates.)

In this Section we analyze more general situation. We consider
relations between  semidensities on an arbitrary odd symplectic
supermanifold and differential forms on even  Lagrangian surfaces
in this \superspace. Then we apply these results for analyzing
relations between conditions (\statement a), (\statement b) and
(\statement c) for Batalin-Vilkovisky formalism geometry.

  Lagrangian surface in
  $(n.n)$-dimensional odd symplectic  supermanifold
  $E=E^{n.n}$
  is $(k.n-k)$-dimensional surface embedded
  in this \superspace such that the
    restriction of symplectic form on it is equal to zero.
    We call $(n.0)$-dimensional Lagrangian surface
    {\it even Lagrangian
    surface.}
    For an odd symplectic \superspace $\PTM$
 an initial underlying $n$-dimensional manifold $M$
  can be considered as
 an even Lagrangian surface
 embedded in this  \superspace.
 (Note that in a case if we consider $\L$-supermanifolds,
 underlying manifold is not necessarily Lagrangian surface.)

   If $L$ is an even Lagrangian surface in odd sympelctic manifold
   $E$ and $\Pi T^*L$ is
   supermanifold associated with cotangent bundle of $L$,
   then one can consider correspondence between
   semidensities on $E$ and differential forms on $L$
   provided there is an identifying symplectomorphism
   between supermanifolds $\Pi T^*L$ and $E$:
\def\identsympl {4.1}
                     $$
\hbox {symplectomorphism}\quad
          \identifyl\colon\quad
   \Pi T^*L\rightarrow E\,\,{\rm and}
                  \quad \identifyl\big\vert_{L}=
                     {\bf id}\,.
                                \eqno (\identsympl)
                          $$
   In this case pull-back $\varphi^*{\bf s}$ of semidensity
    ${\bf s}$ corresponds to differential forms on $L$
    via map (3.4):
   \def\sdensforms {4.2}
                 $$
                   \tau_L^\#\left(\omega_n+w_{n-1}+\dots+\omega_1+\omega_0
                    \right)=
                    \identifyl^*{\bf s}\,,
                            \eqno (\sdensforms) 
                          $$
  where $\omega_k$ is a differential $k$-form on $L$.

 This correspondence depends on a choice of identifying symplectomorphism
 (\identsympl) Thus at first
 we study properties of identifying symplectomorphisms.

\medskip
             {\it 4.1 Identifying symplectomorphisms
             for even Lagrangian surfaces}

             \smallskip

 In usual symplectic calculus
  if $L$ is a Lagrangian surface in a symplectic manifold
   $N$ then there exists symplectomorphism  between
 tubular neighborhoods of $L$ in $T^*L$
 and in $N$ that is identical on $L$ [\Gui].
  In general case there is no
  Lagrangian surface $L$ such that
 $T^*L$ is symplectomorphic to $N$.

 The nilpotency of odd variables leads to the fact that
 odd symplectic supermanifolds have more simple structure.
  Particularly, any $(n.0)$-dimensional surface
  in $(n.n)$-dimensional supermanifold can be expressed locally
  by equations $\theta_i-\Psi_i(x)=0,i=1,\dots,n$
  in any coordinates $\{x^i,\theta_j\}$. Hence
 for every even Lagrangian surface
  in odd symplectic supermanifold $E=E^{n.n}$
  its underlying
  $n$-dimensional manifold
  $M^\prime=M^{\prime n}(L)$
  is an open submanifold in underlying manifold $M$ of $E$.
  If $E^{\prime }$ is a corresponding restriction of
  supermanifold $E$ with underlying manifold $M^{\prime }$
  then one can prove that there exists a symplectomorphism $\varphi$
 that identifies $\Pi T^*L$ with $E^{\prime}$.
  We suppose later that $M^{\prime}$
  coincides with $M$. For example this is a case if $M$
  is a closed connected manifold and $M^{\prime }(L)$
  is also closed. We call such Lagrangian surfaces closed.

     {\bf Proposition 1}
       {\it Let $L$ be an arbitrary closed even
     Lagrangian surface
     in odd symplectic \superspace $E$.
      Then there exists an identifying symplectomorphism}
       (\identsympl)
   {\it between $\Pi T^*L$ and $E$}.

   \smallskip
   Prove this Proposition.

   An identifying symplectomorphism can be
  constructed for every even
  Lagrangian surface in terms of a suitable Darboux coordinates.

   Namely, consider arbitrary atlas
   $\AA(E)=\left[\{x^i_{(\a)},\theta_{j(\a)}\}\right]$
  of Darboux coordinates on a supermanifold $E$
  with closed connected underlying manifold $M$.
 (Every coordinates $\{x^i_{(\a)},\theta_{j(\a)}\}$ of this atlas
 are defined on superdomain $\hat U_{\a}$ with underlying
 domain $U_\a$. Functions $x^i_{0(\a)}$  that are
 numerical parts of functions $x^i_{\a}$,
define an atlas $\left[\{x^i_{0(\a)}\}\right]$ on underlying
manifold $M$.)

 We say that Darboux coordinates $\{x^i,\theta_j\}$ in $E$
 are adjusted to Lagrangian surface $L$
 if $\theta_1=\dots=\theta_n=0$ on $L$.
 Respectively we say that an atlas of Darboux coordinates
  is adjusted to Lagrangian surface $L$ if
  all coordinates from this atlas are adjusted to this surface.

  Suppose that there already exists an identifying symplectomorphism
   $\identifyl$ (\identsympl)
   for a given closed even Lagrangian surface $L$.
   Let  $\AA(\Pi T^*L)=\left[\{y^i_{(\a)},\eta_{j(\a)}\}\right]$
    be an atlas
    of Darboux coordinates in $\Pi T^*L$
   adjusted to cotangent bundle structure of $\Pi T^*L$ (see (3.1)).
   Consider an atlas
   $\AA(E)=\left[\{x^i_{(\a)},\theta_{j(\a)}\}\right]$ of Darboux coordinates
    on $E$ defined by relations
   \def\adjustedatlas {4.3}
               $$
\identifyl^*x^i_{(\a)}=y^i_{(\a)}\,,\quad
\identifyl^*\theta_{j(\a)}=\eta_{j(\a)}\,.
       \eqno (\adjustedatlas) 
              $$
    This atlas is adjusted to
    the Lagrangian surface $L$. Moreover
   from definition of this atlas and (3.1) it follows that
    all transition functions $\Psi_{\a\beta}$ of $\AA(E)$
    on superdomains $\hat U_{\a\beta}$
  with underlying domains $U_{\a\beta}=U_\a\cap U_\beta$
   are "point"-like canonical transformations
   (\cantransform b).
   We also call this atlas on $E$
   {\it an atlas adjusted to
   cotangent bundle structure of Lagrangian surface $L$.}
   It is easy to see that arbitrary atlas of Darboux coordinates
   on $E$ adjusted to cotangent bundle structure of
   a given Lagrangian surface
   $L$  defines some identifying
     symplectomorphism for this Lagrangian surface
     via relations (\adjustedatlas).
    (Darboux coordinates  in r.h.s.
  of (\adjustedatlas)
   adjusted to cotangent bundle structure of $\Pi T^*L$
  are generated by restriction on $L$ of coordinates $\{x^i_{(\a)}\}$
  in $E$.)
  Thus Proposition 1 follows from the following Lemma

   { \bf Lemma 2}. {\it For arbitrary even Lagrangian
   surface $L$ in an odd symplectic supermanifold $E$
    there exists an atlas of Darboux coordinates in $E$
     adjusted to cotangent bundle structure of this surface.}

     Prove this Lemma.

  Considering in a vicinity of arbitrary point of
   $E$ arbitrary Darboux coordinates (see for details Appendix 2)
 we come to some atlas
 $\left[\{x^i_{(\a)},\theta_{j(\a)}\}\right]$
 of Darboux coordinates on $E^{n.n}$.
   If Lagrangian surface $L$ is defined in this atlas by equations
   $\theta_{i(\a)}-\Psi_{i(\a)}(x_{\a})=0$, then
   the condition that surface $L$ is Lagrangian implies that
   $\p_i\Psi_j-\p_i\Psi_j=0$. Hence
   changing $\theta_{i{\a}}\rightarrow$
   $\theta_{i(\a)}-\Psi_{i(\a)}(x_{\a})$ we come to the atlas
   $\AA_{\rm adj}$
   of Darboux coordinates adjusted to the surface $L$
   ($\theta_{i(\a)}\vert_L=0$).

   We show that it is possible to change coordinates in
   every superdomain $U_{\a}$
  for atlas $\AA_{\rm adj}$ in a way that
   all transition functions become "point"-like canonical transformations
   (\cantransform b).
   Prove it by induction.

   Without loss of generality
   consider a case, when a number of charts
    is countable $(\a=1,\dots,n,\dots)$.
   Suppose that we already changed coordinates
  in a required way for first $k$ charts:
 all transition functions $z_{(\a)}=\Psi_{\a\beta}(z_{(\beta)})$
  are already "point"-like canonical transformations
 for $\a,\beta=1,\dots,k$.

 Consider Darboux coordinates
 $\{z^A_{(\a)}\}=\{x^i_{(\a)},\theta_{(\a)}\}$
 on the superdomain  $\hat U_{\a}$
 (with underlying domain
 $U_{\a}$) for $\a=k+1$.
  For every $\beta\leq k$ consider transition function
  (canonical transformation of coordinates)
  $z^A_{(\beta)}=\Psi_{\beta\a}(z_{(\a)})$
  in superdomain $\hat U_{\beta\a}$
  (with underlying domain $U_{\beta\a}= U_{\beta}\cap U_{\a}$).

   All coordinates are adjusted to Lagrangian surface
   ($\theta_{i(\a)}\vert_L=0$), hence from statement 1
   of Lemma 1 it follows
   that one can consider
   in every superdomain $\hat U_{\beta\a}$
   $(\a=k+1,\beta\leq k)$
   new coordinates $\{\t z^A_{\a\beta}\}$ such that
    $z^A_{(\beta)}=\Psi_{\beta\a}(z_{(\a)})=$
    $\F_p\circ\F_{\rm adj}(z_{(\a)})=$
      $z^A_{(\beta)}(\t z_{(\beta\a)}(z_{(\a)}))$
      where $z^A_{(\beta)}(\t z_{(\beta\a)})$
      is point-like canonical transformation and
      $\t z^A_{(\beta\a)}(z_{(\a)}))$
      is adjusted canonical transformation.

  To complete the proof of Lemma we have
  to define in superdomain $\hat U_{\a}$ ($\a=k+1$)
  new coordinates $\{\t z^A_{(\a)}\}$ such that restrictions of
  these coordinates on superdomains $U_{\beta\a}$
  coincide with coordinates $\{\t z^A_{(\beta\a)}\}$
   constructed above. From
   statement 3 of Lemma 1 it follows that there exist
    Hamiltonians $Q_{(\beta\a)}$ in $\hat U_{\beta\a}$
     that generate
    adjusted canonical transformation from coordinates
    $z^A_{\a}$ to coordinates $\t z^A_{(\beta\a)}$
     ($\a\leq k+1,\beta \leq k$).
    From inductive hypothesis and uniqueness of these Hamiltonians
    it follows that
    $Q_{(\a\beta)}=Q_{(\a\gamma)}$ in superdomains
   $\hat U_{\a\beta\gamma}$.
   Hence one can consider
   an odd Hamiltonian obeying to condition $Q=O(\theta^2)$
  on a superdomain $\hat U_{\a}$ ($\a=k+1$) such that restriction
 of this Hamiltonian on superdomains $\hat U_{\beta\a}$
 is equal to $Q_{(\beta\a)}$.
  This Hamiltonian generates
 adjusted canonical transformation
  from coordinates
  $\{z^A_{(\a)}\}$ to a new required Darboux coordinates $\{\t z^A_{(\a)}\}$
  on a superdomain $\hat U_\a$. \finish

\smallskip

\def\E {{\cal E}}
  \def\G{{\cal G}}
  Certainly, the identifying symplectomorphism  (\identsympl) for
   a given closed even Lagrangian surface $L$ is not unique.
   To study this point
 consider (infinite-dimensional) supergroup
 $Can (E)$ of canonical transformations
 of \superspace $E^{n.n}$.
 Every canonical transformation is
  $\L$-point (element) of  this supergroup (see Appendix 1).
  Supergroup $Can(E)$ acts transitively on the
  superspace of closed even Lagrangian surfaces.
  Denote by $Can(L)$ stationary
   subgroup of supergroup $Can(E)$ for $L$
   and consider subgroup $Can_{\rm adj}(L)$
   of supergroup
    $Can(L)$ such that $\L$-points of $Can_{\rm adj}(L)$
  are canonical transformations that are identical on the surface $L$:
  $Can_{\bf adj}(L)\ni F\Leftrightarrow F\vert_L={\bf id}$.
  It is easy to see that canonical transformations
   obeying to this condition
  have following appearance in arbitrary Darboux coordinates adjusted
  to the surface $L$:
  \def\adjincoord {4.4}
                   $$
                   \cases
                   {
   \t x^i=x^i+f^i(x,\theta),
   \quad{\rm where}\,\,f^i(x,\theta)=O(\theta)\,\cr
   \t\theta_i=\theta_i+g_i(x,\theta),\quad
   {\rm where}\,\,g_i(x,\theta)=O(\theta^2)\,\cr
                  }\,,\quad
                   {\rm if}\, \theta_i\vert_L=0\,.
                       \eqno (\adjincoord) 
                    $$
Later we call
   canonical transformations obeying to the condition
   $F\vert_L={\bf id}$
   {\it canonical transformation adjusted to Lagrangian surface $L$}.
    Adjusted canonical transformation
   (\cantransform a) corresponds to transformation
   (\adjincoord) in adjusted coordinates.

 Now consider
  superspace $\Phi(L)$ of identifying symplectomorphisms
   for given closed even Lagrangian surface $L$.
   (Every identifying symplectomorphism $\identifyl$ is $\L$-point
   (element) of this superspace.)
     Supergroup $Can_{\rm adj}(L)$ acts free on superspace
  $\Phi(L)$ of identifying symplectomorphisms:
  arbitrary two identifying symplectomorphisms
  $\identifyl$ and $\identifyl^\prime$ differ on canonical
  transformation adjusted to the surface $L$:
  \def\implies {4.5}
                  $$
  \identifyl^\prime=F\circ\identifyl\,,\quad {\rm where}\,
                          F\Big\vert_L={\bf id}\,.
                         \eqno (\implies)
                          $$

 Consider supergroup $Can_0(E)$ that is unity connectivity component
 of supergroup $Can(E)$, i.e.
 canonical transformation $F$ belongs to $Can_0(E)$
 if it can be included in one-parametric (continuous)
  family $F_t$ of canonical
 transformations ($0\leq t\leq 1$) such that $F_0={\bf id}$
 and $F_1=F$.
   Consider also subgroup $Can_H(E)$ of $Can_0(E)$ such that
   canonical transformation $F$ belongs to $Can_H(E)$
   if it can be included in one-parametric family $F_t$ of canonical
 transformations ($0\leq t\leq 1$)
   generated by some Hamiltonian $Q(x,\theta,t)$:
   $\dot F_t=\{Q,F\}$, $F_0={\bf id}$ and $F_1=F$.
   We call canonical transformations belonging to $Can_H(E)$
   {\it canonical transformations generated by Hamiltonian}.

   Consider Lie superalgebra $\G_{adj}(L)$ such that
     $\Lambda$-points (elements) of this superalgebra
       are odd functions  on $E$
       ("time"-independent Hamiltonians $Q(x,\theta)$)
       that obey to the following condition
       \def\hamdva{4.6}
                          $$
       Q=Q^{ik}(x,\theta)\theta_i\theta_k,\quad
                          {\rm i.e.}\quad
                         Q=O(\theta^2)
                            \eqno (\hamdva) 
                         $$
  in Darboux coordinates  adjusted to Lagrangian surface $L$.
  (Lie algebra structure is defined via odd
  Poisson bracket (\darbouxtheorem).)

One can show that superalgebra $\G_{adj}(L)$ corresponds to
supergroup $Can_{\rm adj}(L)$. Indeed it is  is easy to see that
arbitrary Hamiltonian obeying to condition (\hamdva) generates
one-parametric family  of canonical transformations $F_t=\E xp
\,tQ$ ($0\leq t\leq 1$) adjusted to the surface $L$ and $\E xp
\,tQ_1\not=\E xp \,tQ_2$ if $Q_1\not=Q_2$. (see for detailes
Appendix 3).  Thus the map $\E xp\, Q\colon\,\, \G_{\rm
adj}(L)\rightarrow Can_{\rm adj}(L)$
 is well-defined injection.
 Moreover this exponential map is bijective map.
 To find the Hamiltonian $Q\in \G_{\rm adj}(L)$
    that generates a given
    transformation $F\in Can_{\rm adj}(L)$
    ($F=\E xp\,Q$) consider
      the transformation $F$ in arbitrary atlas $\AA$
    of Darboux coordinates
  adjusted to cotangent bundle structure of the surface $L$.
   In every coordinates from this atlas transformation $F$
   has appearance (\adjincoord), hence according to the statement 3 of Lemma 1
    in every coordinates from atlas $\AA$ there exists unique
     "time"-independent Hamiltonian $Q_{(\a)}$
    obeying to condition (\hamdva) such that
    this Hamiltonian
    generates locally this transformation
    ($Q_{(a)}=-\theta_if^i(x,\theta)+O(\theta^3)$).
    One can see that local Hamiltonians $\{Q_{(\a)}\}$
       do not depend on a choice of coordinates from this atlas.
      Hence they define uniquely a global Hamiltonian $Q$
      in superalgebra $\G_{\rm adj}(L)$. We come to

 {\bf Proposition 2}

  {\it For a given closed even Lagrangian surface $L$ in $E$
  arbitrary two identifying symplectomorphisms
  are related with each other by canonical transformation
  adjusted to Lagrangian surface.
  This canonical transformation
   is generated by "time"-independent Hamiltonian
   that is defined uniquely by condition} (\hamdva).
  {\it In other words supergroup $Can_{\rm adj}(L)$ acts free on superspace
  $\Phi(L)$ of identifying symplectomorphisms.
  The exponential map $\E xp$ from Lie superlalgebra Lie
  $\G_{\rm adj}(L)$ to   $Can_{\rm adj}(L)$ is bijection.}

\smallskip

  For later considerations we need
   to study difference between supergroups
    $Can_0(E)$ (unity connectivity component in $Can(E)$)
    and supergroup $Can_H(E)$ of canonical transformations
    generated by Hamiltonian. For this purpose
  we consider decomposition of
  group $Can(E)$ of all canonical transformations
  on subgroups that are isomorphic to
  $Can_{\rm adj}(L)$, supergroup $Dif\!f (L)$
  of diffeomorphisms of Lagrangian surface $L$,
  and supergroup that acts free
  on the superspace of all even Lagrangian surfaces.

   To describe this latter supergroup
   consider
    abelian supergroup $\Pi Z^1(L)$
   of closed differential one-forms on $L$,
    where $ Z^1(L)$ is superspace of closed
    differential forms on $L$ and
     $\Pi$ is parity reversing functor.
   ($\L$-points of supergroup $\Pi Z^1(L)$ ($Z^1(L)$) are closed one-forms
   with odd (even) coefficients from Grassmann algebra $\Lambda$.)
   Supergroup $\Pi Z^1(L)$ is subgroup of abelian supergroup
   of odd-valued differential one-forms considered in Section 3
   (see (\firstoperation)).
 Superspace $\Pi Z^1(L)$ can be identified with superspace
   of even closed Lagrangian surfaces
   in $\Pi T^*L$, because every
    odd valued differential one-form $\Psi_idx^i$
    can be identified with $(n.0)$-dimensional surface
     embedded
    in $\Pi T^*L$ given by
    equations $\theta_i-\Psi_i(x)=0$.
    Under this identification closed even Lagrangian
    surfaces in $\Pi T^*L$ correspond to closed forms.
    There is a natural monomorphism of supergroup
    $\Pi Z^1(L)$ in supergroup $Can(\Pi T^*(L))$
    of {\it all} canonical transformations of supermanifold $\Pi
    T^*L$: the special canonical transformation (\cantransform c)
      $x^i_{(\a)}\rightarrow x^i_{(\a)}$,
      $\theta_{j(\a)}\rightarrow \theta_{j(\a)}+
      \Psi_{i(\a)}(x_{(\a)})$
    corresponds to an element $\Psi_i(x)dx^i$
      of supergroup $\Pi Z^1(L)$ in an
    atlas of Darboux coordinates on $\Pi T^*L$
    adjusted to cotangent bundle structure of $L$.
    Abelian supergroup $\Pi Z^1(L)$
     acts free on the superspace of
      closed even Lagrangian surfaces in $\Pi T^*L$.
   The action of this supergroup on semidensities in $\Pi T^*L$
  and arbitrary differential forms on $L$ is defined by operation
  (\firstoperation).

    \noindent There is also natural monomorphism of supergroup
    $Dif\!f(L)$ in supergroup
     \noindent $Can(\Pi T^*(L))$
    of {\it all} canonical transformations of supermanifold
    $\Pi T^*L$ corresponding to point-canonical transformation
    (see (\cantransform b) and (3.1)).

\def\AA{{\cal A}}

    Now for supergroups $Dif\!f(L)$ and $\Pi Z^1(L)$
    we consider affine supergroup
    $\Pi Z^1(L)$ $\ltimes$ $Dif\!f(L)$
    such that semidirect product
    is induced by action of diffeomorphisms of $L$
    on forms:
    $[\Psi_1,f_1]\circ[\Psi_2,f_2]=
    [\Psi_1+(f_1^{-1})^*\Psi_2,f_1\circ f_2]$,
    where $\Psi_1,\Psi_2\in\Pi Z^1(L)$ are closed odd valued one-forms and
     $f_1,f_2\in Dif\!f(L)$ are diffeomorphisms of $L$.
    Monomorphisms  of supergroups $\Pi Z^1(L)$ and $Dif\!f(L)$
    in supergroup $Can(\Pi T^*L)$ considered above
     define monomorphism $\iota$
     of the affine supergroup
     $\Pi Z^1(L)$ $\ltimes$ $Dif\!f(L)$ in
    the supergroup $Can(\Pi T^*L)$.
     Thus every identifying symplectomorphism $\identifyl$
     defines monomorphism
     $\iota_{\identifyl}=$ $\identifyl\circ\iota\circ\identifyl^{-1}$
     of the supergroup $\Pi Z^1(L)$ $\ltimes$ $Dif\!f(L)$
     in the supergroup $Can(E)$.

On the other hand consider for arbitrary canonical transformation
$F\in Can(E)$ Lagrangian surface $\widetilde
L=\identifyl^{-1}\circ F(L)$ in $\Pi T^*L$
 and closed odd valued one-form $\Psi$
 on $L$ corresponding to the Lagrangian surface $\widetilde L$.
 Then canonical transformation $F^\prime=$
 $\iota_{\identifyl}([-\Psi,{\bf id}])\circ F$
 of \superspace $E$ belongs to supergroup
$Can(L)$ of canonical transformations that transform Lagrangian
surface $L$ to itself. The restriction of the canonical
transformation $F^\prime$ on $L$ defines diffeomorphism
$f=F^\prime\vert_L\in Dif\!f(L)$. Thus we define projection map
\def\projectionmap{4.7}
                    $$
p_{\identifyl}\colon\quad Can(E)\rightarrow \Pi
Z^1(L)\hbox{$\ltimes$} Dif\!f(L)\,
                  \eqno(\projectionmap)
                  $$
 that depends on identifying symplectomorphism.

This map projects subgroup $Can_{\rm adj}(L)$ of canonical
transformations adjusted to the surface $L$ to unity element and
obeys to condition $p_{\identifyl}\circ\iota_{\identifyl}={\bf
id}$.
 We come to the following result:
 for a given Lagrangian surface $L$ and identifying symplectomorphism
 $\identifyl$ arbitrary canonical transformation $F\in Can(E)$
 can be decomposed uniquely in the following way:
   \def\newdva{4.8}
                         $$
                F=\iota_{\identifyl}([\Psi,f])\circ F_{\rm adj}=
              F_s\circ F_p\circ F_{\rm adj}\,,\quad{\rm where}
                          $$
                          $$
      [\Psi,f]=p_{\identifyl}(F)\,,\,\,
            F_s=\iota_{\identifyl}([\Psi,{\bf id}])\,,\,\,
             F_p=\iota_{\identifyl}([0,f])\,,\,\,
             F_{\rm adj}\in Can_{\rm adj}(L)\,.
              \eqno (\newdva) 
                          $$
(The decomposition (\cantransform) in Lemma 1 corresponds to this
decomposition.)

 One can check the following property of projection map (\projectionmap):
  If $\identifyl$, $\identifyl^\prime$ are two arbitrary
  identifying symplectomorphisms for a given Lagrangian surface
  $L$ and $p_{\identifyl}(F)=[\Psi,f]$,
  $p_{\identifyl^\prime}(F)=[\Psi^\prime,f^\prime]$ then
  \def\projectionchanging {4.9}
                       $$
     \Psi^\prime-\Psi=d\Phi\,, \quad f^\prime=f_0\circ f\,,
                      \eqno (\projectionchanging)
                      $$
  where $f_0$ is diffeomorphism that can be included in one-parametric
  continuous family $f_t$
   of diffeomorphisms such that $f_1=f$ and $f_0={\bf id}$.
   In other words $f_0$ belongs to group $Dif\!f_0(L)$ that is unity
   connectivity component of group $Dif\!f(L)$.

   According to Proposition 2 conditions (\projectionchanging)
  have to be checked only for infinitesimal canonical transformations
  (\adjincoord) adjusted to Lagrangian surface $L$ and
  generated by Hamiltonian (\hamdva).
  This can be done by easy straightforward calculations.

 From (\projectionchanging) it follows that
  for a given Lagrangian surface $L$
projection map (\projectionmap) defines a map
\def\projectionmaphom{4.10}
                   $$
   p_L\colon\quad Can(E)\rightarrow
   \Pi H^1(L)\hbox{$\ltimes$}\pi_0(Dif\!f(L))\,,
                  \eqno(\projectionmaphom a) 
               $$
    where $\Pi H^1(L)$ is abelian group of cohomology classes
     of one-forms on $L$ (with reversed parity)
     and  $\pi_0(Dif\!f(L))=Dif\!f(L)/Dif\!f_0(L)$
     is discrete group of connectivity components of $Dif\!f(L)$.
 Using decomposition (\newdva),
  relations (\projectionchanging)
  and the fact that supergroup $Can_{\rm adj}$ is normal subgroup in $Can(L)$
  one can show that (\projectionmaphom a) is epimorphism.
 (Projection map (\projectionmap) is not epimorphism, because
 supergroup $Can_{\rm adj}$ is not normal subgroup in $Can(E)$.)

 One can consider also a composition of epimorphism (\projectionmaphom a)
 with natural epimorphism of $\Pi H^1(L)$ $\ltimes$ $\pi_0(Dif\!f(L))$
  on  $\pi_0(Dif\!f(L))$:
                       $$
                       \hat p_L\colon\quad
  Can(E){\buildrel p_L\over\longrightarrow}
   \Pi H^1(L)\hbox{$\ltimes$}\pi_0(Dif\!f(L))\rightarrow\pi_0(Dif\!f(L))\,.
                  \eqno(\projectionmaphom b) 
                  $$
  Epimorphisms (\projectionmaphom a) and (\projectionmaphom b)
  allow to check difference between supergroups
 $Can_0(E)$ and $Can_H(E)$ because
 \def\kernels {4.11}
                     $$
            {\bf ker}\, p_L=Can_H(E)\,\quad {\rm and}\quad
            {\bf ker}\,\hat p_L=Can_0(E)\,.
                          \eqno (\kernels)
                    $$
 Namely consider arbitrary canonical transformation
 $F$ that belongs to the kernel of epimorphism (\projectionmaphom a).
 Then for projection map (\projectionmap) $p_{\identifyl}(F)$
 $=[\Psi,f]$ where $\Psi=d\Phi$ and $f\in Dif\!f_0(L)$.
 Consider decomposition (\newdva) for this canonical transformation $F$.
 Then
 canonical transformation $F_s=\rho([\Psi,{\bf id}])$
 is generated by Hamiltonian $Q=\Phi(x)$.
 Canonical transformation
 $F_s=\rho([0,f])$ is generated by Hamiltonian $Q=K^i(t,x)\theta_i$
 where "time"-dependent vector field
   $K^i(t,x)$ is equal to $f^{-1}_t\circ\dot f_t$
   for a family $f_t$ of diffeomorphisms that
   connects diffeomorphism $f$ with identity diffeomorphism.
  Canonical transformation $F_{\rm adj}$ is generated by some Hamiltonian
  $Q(x,\theta)$
 according to Proposition 2.
 Hence the kernel of epimorphism
 (\projectionmaphom a) belongs to $Can_H(E)$.
   To prove the converse implication consider one-parametric family
   $F_t$ of canonical transformations generated by arbitrary Hamiltonian
   $Q(x,t)$. Decompose for every $t$ transformation $F_t$
   by formula (\newdva) for arbitrary identifying
   symplectomorphism $\identifyl$:
   $F_t=$ $F_s(t)\circ F_p(t)\circ F_{\rm adj}(t)$.
   Transformations $F_p(t)$ and $F_{\rm adj}(t)$ are generated by
     Hamiltonians hence transformation $F_s(t)$ is generated
     by Hamiltonian $Q^\prime$ also. Hence $\Psi=dQ^\prime$,where
     $p_{\identifyl}F_p(t)=[\Psi_t,1]$ and
     $p_L(F)=[0,1]$ in $\Pi H^1(L)$ $\ltimes$ $\pi_0(Dif\!f(L))$.

   The proof of the second relation in (\kernels) is analogous.
We come to

{\bf Proposition 3} {\it Let $L$ be a closed even Lagrangian
surface in an odd symplectic supermanifold $E$. Let $Can_0(E)$ be
unity connectivity component of supergroup
 $Can(E)$ of canonical transformations of $E$ and $Can_H(E)$ be
  supergroup of canonical transformations generated by Hamiltonian.
   Then the following relations
 between supergroups $Can(E)$, $Can_0(E)$ and $Can_H(E)$ are obeyed:
                         $$
                         \matrix
                         {
  Can(E)/Can_0(E)=\pi_0(Dif\!f (L)),\cr
 Can(E)/Can_H(E)=\Pi H^1(L)\hbox{$\ltimes$}\pi_0(Dif\!f(L)),\cr
         Can_0(E)/Can_H(E)=\Pi H^1(L)\,.\cr
                       }
                    $$
In particularly supergroup $Can_0(E)$ is equal to supergroup
 $Can_H(E)$ if $H^1(L)=0$.}

\smallskip

    \noindent Groups $Dif\!f(L)$, $\pi_0(Dif\!f(L)$), $\Pi Z^1(L)$ and $\Pi H(L)$
      are isomorphic to groups $Dif\!f(M)$, $\pi_0(Dif\!f(M)$),
       $\Pi Z^1(M)$ and $\Pi H(M)$
      respectively,
      where $M$ is underlying supermanifold,
      but isomorphisms
      are not canonical.

\smallskip
             {\it 4.2 Relation between semidensities and differential forms
              on a Lagrangian surface}

             \smallskip

    Now we return to relation (\sdensforms)
    between semidensities on odd symplectic supermanifold
    $E=E^{n.n}$ and differential forms on even Lagrangian surfaces.
   We assume that underlying manifold is orientable
    (see Remark after (\bercanonical)) and its orientation
    is fixed.
    This fixes orientation on even Lagrangian surfaces.

        We note also that if we consider the points of supermanifold as
  $\L$-points  where $\L$ is an arbitrary Grassmann algebra,
  then one have to consider differential forms with coefficients
   in this algebra $\L$ (see Appendix 2). It follows from
   (3.4) that if $\bf s$
   is even (odd) semidensity then $k$-form in l.h.s. of relation
   (\sdensforms) has coefficients in Grassmann algebra $\L$ with
    parity $p=(-1)^{n-k}$
   ($p=(-1)^{n-k+1}$).
   More precisely denote by $S$ a superspace of
   semidensities in $E$.
  $\L$-points of superspace $S$ are even semidensities, i.e.
  semidensities ${\bf s}=s(x,\theta)\sqrt{D(x,\theta)}$,
   such that
   $s(x,\theta)$ are even functions with coefficients in
   Grassmann algebra $\L$.
  Denote by $\Omega^k$
   superspace of differential $k$-forms on even Lagrangian surface
   $L$ and consider also superspace $\Pi\O^k$,
   where $\Pi$ is parity reversing functor.
   $\L$-points of superspace $\O^k$ are
  differential $k$-forms
  with even coefficients
   from Grassmann algebra $\L$,
   $\L$-points of superspace $\Pi\O^k$ are differential $k$-forms,
  with odd coefficients from Grassmann algebra $\L$.
   Consider a superspace
\def\spaceV{4.12}
   $$
   \O^*(L)=\O^n\oplus\Pi \O^{n-1}\oplus\O^{n-2}\oplus\Pi\O^{n-3}
   \oplus\O^{n-4}\dots
    \eqno (\spaceV) 
      $$
         Relation (\sdensforms) defines a map:
    \def\mapA {4.13}
                   $$
      w(L,\identifyl,{\bf s})=(\tau_L^{\#})^{-1}\identifyl^*{\bf s}\,
      \eqno (\mapA) 
                    $$
         between superspace $S$
       and superspace $\O^*(L)$.
       (Here and later where it will not lead to confusion
       we denote by $w$ a linear combination of differential
       forms $w_{n}+w_{n-1}+\dots+w_{0}$.)
     At what extent map (\mapA)
  depends on a choice of identifying symplectomorphism
   and on a choice of even Lagrangian surface?

    If $F$ is arbitrary canonical transformation of
   $E$ and $\identifyl^\prime=F\circ\identifyl$
    then for map (\mapA)
   \def\propertyofmapA{4.14}
                        $$
   w(L,\identifyl^\prime,{\bf s})=
            w(L,F\circ\identifyl,{\bf s})=w(L,\identifyl,F^*{\bf s})\,.
             \eqno (\propertyofmapA)
                      $$
    Thus bearing in mind Proposition 2
    we study  the action of
    supergroup $Can_{\rm adj}(L)$ of canonical transformations
    on semidensities.

\smallskip

  {\bf Proposition 4}

   a){\it Let ${\bf s}$ be arbitrary semidensity on
   odd symplectic supermanifold
   $E=E^{n.n}$ with closed connected underlying manifold $M^n$
   and  $F$ be arbitrary canonical transformation of $E^{n.n}$
   adjusted to a given
   even Lagrangian surface $L$ in $E$ ($F\big\vert_L={\bf id}$,
    i.e.$F\in Can_{\rm adj}(L)$).

   Then} $\left(F^*{\bf s}-{\bf s}\right)\big\vert_L=0$.

\noindent b) {\it Arbitrary canonical transformation $F$
 generated by Hamiltonian
    ($F\in Can_H(E))$ chan\-ges arbitrary closed
    semidensity on an exact form:
   if $\Ds {\bf s}=0$ then $F^*{\bf s}-{\bf s}=\Ds{\bf r}$ .
   In the case if this transformation is adjusted to Lagrangian surface $L$
   ($F\in Can_{\rm adj}(L)$ $\subseteq Can_H(L)$), then
   condition $\left(F^*{\bf s}-{\bf s}\right)\big\vert_L=0$ is obeyed also.}

\noindent c) {\it If ${\bf s}$ and ${\bf s}_1$ are
 arbitrary even closed non-degenerate semidensities (
   ${\bf s}, {\bf s}_1\in \B_{\rm deg}$ (see} (\bvsdensdef)),
  {\it differ on an exact semidensity: ${\bf s}_1-{\bf s}=\Ds {\bf r}$,
   and for even Lagrangian surface $L$ condition
    $({\bf s}_1-{\bf s})\big\vert_L=0$ is obeyed,
   then there exists a canonical transformation $F$
   adjusted to $L$ such that ${\bf s}=F^*{\bf s}_1$.}

\smallskip

(We say that semidensity ${\bf s}$ is equal to zero on even
Lagrangian surface $L$ (${\bf s}\big\vert_L=0$)
 if in Darboux coordinates adjusted to $L$
${\bf s}$ $=s(x,\theta)\sqrt{D(x,\theta)}$ with
$s(x,\theta)\vert_{\theta=0}=0$.)

 The statement a) follows
  from explicit expression (\adjincoord) for transformation $F$
  adjusted to Lagrangian surface $L$.
   According to Proposition 2 statement b) have to be checked
    only for infinitesimal transformations generated by Hamiltonian.
    For these transformations
    this statement follows from formula (\denstransform).

   To prove statement c) we consider
   \def\hamtimedepend {4.15}
   a following "time"-depending Hamiltonian:
                     $$
             Q(t)={-{\bf r}\over{\bf s}+t\Ds{\bf r}}\,,\quad 0\leq t\leq 1
                           \eqno (\hamtimedepend)
                        $$
          for any one-parameter family
          ${\bf s}_t={\bf s}_0+t\Ds{\bf r}$, $0\leq t\leq 1$
          of even closed non-degenerate semidensities
          (${\bf s}_t\in \B_{\rm deg}$
          at any $t$).

    It is easy to check that during a "time" $t$
    this Hamiltonian generates canonical transformation
    $F_t$ that transforms ${\bf s}_t$ to ${\bf s}$
     ($F_t^*{\bf s}_t={\bf s}$).
     Indeed according  to (\hamtimedepend) and (\denstransform)
    if transformation $F_t$ obeys to conditions
    $\dot F_t=\{Q,F_t\}$ and $F_0$ $={\bf id}$ then
                          $$
    {d\over dt}F^*_t{\bf s_t}=F^*_t\Ds{\bf r}+
                      F_t^*\left(\Ds\left( Q(t){\bf s}_t\right)\right)=0\,
                      \Rightarrow F_t^*{\bf s}_t={\bf s}_0\,.
                            $$

\def\G{{\cal G}}
 Consider Hamiltonian (\hamtimedepend) for semidensities
  ${\bf s}, {\bf s}_1\in\B_{\rm deg}$  with
 $\Ds{\bf r}={\bf s}_1-{\bf s}$ choosing ${\bf r}$
  in such a way that ${\bf r}=O(\theta^2)$ in coordinates
  adjusted to $L$.
 Then Hamiltonian (\hamtimedepend) leads
  to canonical transformation
  $F_t$ that is adjusted to $L$ at any $t$
  and the transformation
  $F=F_1$ transforms ${\bf s}_1$ to ${\bf s}$.  \finish

\smallskip

Now we use this Proposition for analyzing relation (\sdensforms)
for a given even ($(n.0)$-dimensional) Lagrangian surface $L$.

1. According to Proposition 2
  two identifying symplectomorphisms for a given even Lagrangian surface
 differ on canonical transformation
 adjusted to this surface. Hence from statement a) of Proposition 4
and condition (\propertyofmapA) for map (\mapA) it follows that
top-degree form $w_n$ in (\mapA) does not depend on a choice of
identifying symplectomorphism $\identifyl$:
 for a given even Lagrangian surface $L$ relations
(\mapA) induce a well-defined map
\def\volumemap{4.16}
              $$
V(L, {\bf s})=w_n(L, {\bf s}){\buildrel\,{\rm def}\over=}\,
w_n(L,\identifyl,{\bf s})\,,
                  \eqno(\volumemap)
              $$
 where $\identifyl$ is arbitrary identifying symplectomorphism $\identifyl$
 for Lagrangian surface $L$.

 Formula (\volumemap) defines the
  map from superspace $S$ of semidensities in $E=E^{n.n}$
 to the superspace
 $\Omega^n(L)$ of top-degree forms on $L$.
 This means that
 semidensity can be considered as
 well-defined integration object over even
 Lagrangian surface.
This corresponds to general result that
 semidensity can be
considered as
 an integration object over arbitrary
 $(n-k.k)$-dimensional Lagrangian surface.
 (See [\SchCMPa] for corresponding construction and
 [\JMPd], [\tyutinb] for explicit formulae.)

2. Consider restriction of the map (\mapA) on the superspace
  $\B$ of closed semidensities. If
 semidensity ${\bf s}$ is closed then it
 follows from statement b)
of Proposition 4 and relation (\innerdelta) that  for a given
even Lagrangian surface $L$ under changing of identifying
symplectomorphism
  $\varphi_{_{L,E}}\mapsto$
  $\varphi^\prime_{_{L,E}}=F_{\rm adj}\circ\varphi_{_{L,E}}$,
  corresponding differential forms in (\mapA)
   change on exact forms: if $\Ds{\bf s}=0$ then
\def\changingdva {4.17}
                  $$
  w_k(L,\identifyl^\prime,{\bf s})=w_k(L,F\circ\identifyl,{\bf s})=
        w_k(L,\identifyl,F^*{\bf s})=w_k(\identifyl,{\bf s})+
                       d w_{k-1}(L,\identifyl,{\bf r})\,.
              \eqno (\changingdva) 
                       $$
In particular  $w_0$ (constant) as well as $w_n$ do not depend on
identifying symplectomorphism.

Projection of  superspaces $\O^k$ in (\spaceV) on a superspaces
$H^k$ of cohomology classes for $k\leq n-1$ induces projection of
superspace $\O^*$ on superspace:
\def\spaceVcoh{4.18}
                     $$
       \O^n(L)\oplus
       \Pi H^{n-1}(L)\oplus
           H^{n-2}(L)\oplus
       \Pi H^{n-3}(L)\oplus
           H^{n-4}(L)\dots
    \eqno (\spaceVcoh) 
                        $$
 From (\changingdva) it follows that
  considering map (\mapA) on closed semidensities
  for arbitrary identifying symplectomorphisms
 and projecting value of this map on the superspace (\spaceVcoh)
  we come to well-defined  map
\def\mapcoh {4.19}
                 $$
\hat V(L,{\bf s})=
 w_n(L,{\bf s})+[w_{n-1}](L,{\bf s})+\dots+[w_0](L,{\bf s}),
 \quad{\rm if}
              \,\,\Ds{\bf s}=0
 \eqno (\mapcoh) 
                    $$
($[w_k]$ is cohomology class of form $w_k$).
 The map $\hat V(L,{\bf s})$ is linear surjection map
  from the space
  $\S$ of closed semidensity on superspace (\spaceVcoh).

 By definition $\hat V(L,{\bf s})$ is $Can_{\rm adj}$-invariant map:
it does not change under arbitrary canonical transformation
adjusted to the surface $L$:
\def\bijectivity {4.20}
                         $$
                    {\rm if}\,\,
 {\bf s}_1=F^*{\bf s}_2\,\,{\rm where}\,\,
 F\in Can_{\rm adj}(L)\,\, {\rm then}\,\,
 \hat V(L, {\bf s}_1)=\hat V(L,{\bf s}_2)\,.
                                \eqno(\bijectivity)
                           $$
   The opposite implication is obeyed
  for closed non-degenerated semidensities.
  Namely consider map (\mapcoh) for the
  subset $\B_{\rm deg}$ of closed even non-degenerate semidensities
  (see \bvsdensdef).
   If $\hat V(L,{\bf s}_1)=\hat V(L,{\bf s}_2)$
   for two arbitrary closed non-degenerate semidensities, then
  from statement c) of Proposition 4
  it follows that there exists canonical transformation
  $F$ adjusted to surface $L$ such that
  ${\bf s}_1=F^*{\bf s}_2$.

\smallskip

  Now we analyze dependence of
  map  (\mapcoh) under a changing of even Lagrangian surface $L$.
 We study this point from more general point of view
 considering an action of group $Can(E)$ of all
 canonical transformations on maps (\mapA) and (\mapcoh).

Projecting superspace $\O^n(L)$ of top-degree forms on superspace
$H^n(L)$ of corresponding cohomology classes we come from the map
$\hat V(L,{\bf s})$ to the map
\def\cohomologymap {4.21}
               $$
               \hat H(L,{\bf s})=
               [w](L,{\bf s})=
               [w_n](L,{\bf s})+\dots+[w_0](L,{\bf s})
                  \eqno (\cohomologymap) 
               $$
that is defined on superspace $\B$ of closed semidensities and
takes values in a superspace
                      $$
      H^*(L)=H^n(L)\oplus\Pi H^{n-1}(L)\oplus
      H^{n-2}(L)\oplus\Pi H^{n-3}(L)\oplus\dots
                     $$
From statement b) of Proposition 4 it follows that this map is
$Can_H(E)$-invariant: it does not change under arbitrary
canonical transformations generated by Hamiltonian:
\def\cohmapinvariant {4.22}
                      $$
{\rm if}\,\, {\bf s}_1=F^*{\bf s}_2\,\,{\rm where}\,\, F\in
Can_H(E)\,\,{\rm then}\,\,
        \hat H(L,F^*{\bf s})=\hat H(L,{\bf s})\,.
                           \eqno(\cohmapinvariant) 
                       $$
The opposite implication is obeyed under the following
restriction.
 Let ${\bf s}_1$ and ${\bf s}_2$ be
   two closed non-degenerated semidensities
   (${\bf s}_1,{\bf s}_2\in \B_{\rm deg}$) such that
   for the map (\cohmapinvariant)
   $\hat H(L,{\bf s}_1)=\hat H(L,{\bf s}_2)$.
   In this case ${\bf s}_2={\bf s}_1+\Ds {\bf r}$.
   If one-parametric family of closed semidensities
   ${\bf s}_t={\bf s}_1+t\Ds {\bf r}$ $(0\leq t\leq 1)$
   belongs also to $\B_{\rm deg}$ then
   there exists canonical
transformation $F$ generated by Hamiltonian ($F\in Can_H(E))$
such that $F^*{\bf s}_1={\bf s}_2$. One comes to this
transformation considering Hamiltonian (\hamtimedepend)\footnote
{$^*$}
  {This leads to the statement
  of Theorem 5 in the paper [\SchCMPa]
  in the special case when
  $Can_0(E)=Can_H(E)$, i.e. $H^1(M)=0$ (see Proposition 3).}.

  Consider now the action of arbitrary canonical transformation
  on map (\cohomologymap).

 From Proposition 3 and (\cohmapinvariant)
 it follows that
 under arbitrary canonical transformation the map
 $\hat H(L.{\bf s})$
 have to be transformed under the action of the group
 $Can(E)/Can_H(E)=$
 $\Pi H^1(L)$$\ltimes$$ \pi_0(Dif\!f(L))$.
 Namely
 \def\svjazjform {4.23}
                      $$
 \hat H(L,F^*{\bf s})= [f]^*
 \left([\Psi]\,{\cal d}\,\hat H(L,{\bf s})\right)\,,
                              \eqno(\svjazjform)
                       $$
where  $\left[[\Psi],[f]\right]$ is an element of supergroup $\Pi
H^1(L)$$\ltimes$$ \pi_0(Dif\!f(L))$ defined by the action of
epimorphism (\projectionmaphom a) on canonical transformation $F$
  and
the operation $\,{\cal d}\,$
 is defined for semidensities and corresponding differential
 forms by operations (\firstoperation) and (\firstoperationa).
The pull-back  $[f]^*$ of equivalence class $[f]$ is well-defined,
because pull-back $f^*_0$ of
 diffeomorphism $f_0\in Dif\!f_0(L)$ acts identically
 on cohomological classes of differential forms.

On the other hand it follows from (\propertyofmapA) that
 \def\anothersurface{4.24}
                          $$
\hat H(L,F^*{\bf s})= (F\vert_L)^*\hat H(\widetilde L,{\bf s})\,,
                            \eqno (\anothersurface)
                $$
where $\widetilde L$ is an image of Lagrangian surface $L$ under
canonical transformation $F$,

  One can easy derive formulae (\svjazjform) and (\anothersurface)
 from (\mapA)
 performing calculations in arbitrary Darboux coordinates adjusted
 to cotangent bundle structure of Lagrangian surface $L$
 (i.e. choosing arbitrary identifying symplectomorphism $\identifyl$)
 and using decomposition formula (\newdva).

  It is useful to rewrite formulae (\svjazjform) and (\anothersurface)
 in components:
 \def\svjazjformcom {4.25}
                       $$
[w_k](L,F^*{\bf s})=(F\vert_L)^*[w_k](\widetilde L,{\bf s})=
            [f]^*\left(
            \sum_{p=0}^k {1\over p!}
    \underbrace {[\Psi]\wedge\dots\wedge [\Psi]}_{p\quad{\rm times}}\wedge
                          [w_{k-p}]
                          \right)
                             \,.
                                \eqno (\svjazjformcom)
                   $$

We note that if for a given pair $(L,\widetilde L)$ of even
Lagrangian surfaces canonical transformation $F$ transforms $L$
to $\widetilde L$
  then cohomological class of odd valued one-form
  corresponding to the pair $(F,L)$ by epimorphism
  (\projectionmaphom a) is well-defined function of the pair
   $(L,\widetilde L)$:
 \def\twosurfacefunction {4.26}
                       $$
             \Pi H^1(L)\ni[\Psi]=[\Psi](L,\widetilde L)\,.
                               \eqno(\twosurfacefunction)
                        $$
The pull-back $(F\vert_L)^*$ of restriction $F\vert_L$ of
canonical transformation $F$ on Lagrangian surface $L$ induces
bijective map between differential forms and corresponding
cohomological classes on surfaces $L$ and $\widetilde L$. Using
(\svjazjformcom) we can compare cohomological classes of
differential forms corresponding to a given closed semidensity for
two different even Lagrangian surfaces.
 In particular, from  (\svjazjformcom)
  it follows that for arbitrary closed semidensity ${\bf s}$
and for arbitrary pair of closed Lagrangian surfaces
$(L,\widetilde L)$
\def\svjazjforma {4.27}
                   $$
        [w_k](\widetilde L,{\bf s})=0\quad{\rm if}\quad
          [w_0]( L,{\bf s})=\dots=[w_k](L,{\bf s})=0\,,
                      $$
and in the case if cohomological class of one-form
 $[\Psi](L,\tilde L)$ in (\twosurfacefunction) is equal to zero, then
           $$
        [w_k](\widetilde L,{\bf s})=0\quad{\rm iff}\quad
          [w_k]( L,{\bf s})=0\,.
            \eqno ( \svjazjforma)
              $$

  The simple but important consequence of
   these considerations is following:
     $[w_0]$-component of function $\hat H(L,{\bf s})$
  (\cohomologymap) does not depend on canonical transformation
 and it is invariant constant on all Lagrangian surfaces.

   {\bf Corollary 1}
   {\it To every closed semidensity ${\bf s}$ ($\Ds {\bf s}=0$)
    corresponds a positive constant
   $c({\bf s})$.
     If in arbitrary Darboux coordinates}
                     $$
  {\bf s}=s(x,\theta)\sqrt{D(x,\theta)}=
  (\rho(x)+b^i(x)\theta_i+\dots+
         c\theta_1\theta_2\dots\theta_n)\sqrt{D(x,\theta)}
                        $$

 {\it then $c({\bf s})=|c|$.
 This constant does not depend on the choice of Darboux coordinates
 and on the changing of density under arbitrary canonical transformation.
 This constant is equal (up to a sign)
        to cohomological class  $[w_0]$
      of zeroth order differential form
      corresponding to semidensity ${\bf s}$
      on arbitrary even Lagrangian surface $L$.
           (A sign of $c({\bf s})$ depends on orientation.)}

\smallskip

   Note that
$c({\bf s})$ can be considered as integral of semidensity
 ${\bf s}$ over Lagrangian $(0.n)$-dimensional surface
 $x^1=x^1_0,\dots, x^n=x^n_0$:
  $c=\int_L{\bf s}$ $=\int s(x_0,\theta)d^n\theta$.

\smallskip

             {\it 4.3 Application to BV-geometry}

             \smallskip

Now using results obtained in this Section we analyze Statement 1
(see Introduction)
 of Batalin-Vilkovisky master equation.

  Let ${\bf s}$ be an arbitrary closed semidensity in $\B_{\rm deg}$,
  i.e. non-degenerate semidensity that obeys to BV-master equation
   (\statementb, \bvsdensdef).

   For arbitrary $\L$-point $\a$ in the odd symplectic supermanifold
   $E=E^{n.n}$ consider
   arbitrary closed even Lagrangian surface $L$
   such that this point belongs to this surface
   and choose arbitrary identifying symplectomorphism $\identifyl$,
   corresponding to this surface,
   i.e. atlas of Darboux coordinates adjusted to cotangent bundle structure
   of this Lagrangian surface.
   Consider on $L$
   differential form
\def\bvodin{4.28}
           $$
   w(L,\identifyl,{\bf s})=w_n+w_{n-1}+\dots+w_0
   \eqno (\bvodin)
   $$
   defined by the map (\mapA).

  Locally all closed differential
  forms except zeroth-forms are exact
  and $[w_0]=\pm c({\bf s})$ is invariant constant according
  to Corollary 1.
  Hence using statement c) of Proposition 4
  one can find canonical transformation adjusted to this Lagrangian surface
  and correspondingly another identifying symplectomorphism
      $\identifyl^\prime$
  such that in (\bvodin) all differential forms
  $w_k$ for $k=1,\dots,n-1$
  vanish in a vicinity of the point $\a$.
  Consider Darboux  coordinates \Darbouxcoord
  on $E$ in a vicinity of this point from
  the atlas of Darboux coordinates corresponding to the
   identifying symplectomorphism  $\identifyl^\prime$ and
   by suitable "point" canonical transformation
   (\cantransform b) choose them in a way
   that differential form $w_n$ is equal to $dx^1\wedge\dots\wedge dx^n$
   in these coordinates. Thus we come to Darboux coordinates
   in a vicinity of the point $\a$ such that in these Darboux coordinates
 semidensity ${\bf s}$
has following appearance:
\def\bvsdensnormalform {4.29}
                    $$
  {\bf s}=s(x,\theta)\sqrt{D(x,\theta)}=
  (1+c\theta_1\theta_2\dots\theta_n)\sqrt{D(x,\theta)}\,,
                  \eqno(\bvsdensnormalform)
                    $$
   where $c$ is equal  up to sign to the invariant
    constant $c(\bf s)$ corresponding to the semidensity ${\bf s}$
    (see Corollary 1).
  The condition
   $c({\bf s})\not=0$ is the obstacle to condition (\statementa).

     Consider now the value of the map (\mapcoh)
     on this Lagrangian surface:
\def\bvdva{4.30}
                          $$
                  \hat V(L,{\bf s})=
                  w_n+[w_{n-1}]+\dots+[w_0]\,.
                                    \eqno (\bvdva)
                         $$
  If $\hat V(L,{\bf s)}=w_n+c$, i.e. all cohomological classes
   $[w_k]$ for $k=1,\dots,n-1$ in (\bvdva) vanish on the surface $L$,
    then
  one can consider identifying symplectomorphism $\identifyl$
   such that
  $\tau^\#_L(w_n+c)=\identifyl^*{\bf s}$ for the map (\mapA).
    It means that there exists
   an atlas $\left[\{x^i_{(\a)},\theta_{j(\a)}\}\right]$
   of Darboux coordinates on $E^{n.n}$ adjusted to cotangent
   bundle structure
    of Lagrangian surface $L$
   such that
   in arbitrary coordinates from this atlas semidensity
    ${\bf s}$ is expressed by relation (\bvsdensnormalform).
   Semidensity ${\bf s}$ has appearance
    $\sqrt{D(x,\theta)}$
    in any Darboux coordinates from this atlas
    if invariant constant $c({\bf s})=0$.
    In other words in this case supermanifold can be identified with
    $\Pi T^*L$ with volume form on $\Pi T^*L$ induced by
    volume form on $L$.

     It follows from (\svjazjform ---\svjazjforma) that
     this statement holds for another even Lagrangian surface
     $\t L$ iff cohomological class $[{\Psi}]$ of
     one-form corresponding
      to a pair $(L,\t L)$ of Lagrangian surfaces
      (see \twosurfacefunction)
       is equal to zero. In particular
     this statement is irrelevant to a choice of Lagrangian
     surface if $H^1(M)=0$. ($M$ is underlying manifold for $E$.)

  Now we analyze condition (\statementc)
  for even-nondegenerate semidensity ${\bf s}=\sqrt{d\bf v}$.
  From (\deltadenspropb ) it follows that condition
     (\statementc) means that function (\denscalculus d)
     is equal to an odd constant $\nu$, and
     $\Ds {\bf s}=\nu{\bf s}$.
     One can see using correspondence between semidensities
     and differential forms that
     all solutions to this equation
    are  following: ${\bf s}=\Ds {\bf h}-\nu {\bf h}$,
   where ${\bf h}$ is an arbitrary semidensity.
The odd constant $\nu\not=0$
     is the obstacle to condition (\statement b),
  if condition (\statementc) is obeyed.

 \noindent We come to the

   {\bf Corollary 2}

   {\it Let $E=E^{n.n}$
   be an odd symplectic supermanifold with connected orientable
   underlying manifold $M$ and this supermanifold is
   provided with a volume form
   $\dv$, such that $\D^2=0$. Then}

   1) {\it to the volume form $\dv$ corresponds the odd constant
   $\nu$: $\Ds\sqrt\dv=\nu\sqrt\dv$ and
   $\sqrt\dv$
   $=\Ds {\bf h}-\nu {\bf h}$}
   for some odd semidensity ${\bf h}$.

   2) {\it If the odd constant $\nu$ is equal to zero, then
   the master-equation $\Ds\sqrt\dv=0$ holds for semidensity $\sqrt\dv$.
   In this case to the volume form $\dv$ corresponds non-negative constant
   $c=c(\sqrt\dv)$ and there exists an atlas of Darboux coordinates
   on $E^{n.n}$ such that $\dv=(1\pm 2c)D(x,\theta)$
    in any coordinates
   from this atlas.}

   3) {\it In the case if the constant $c({\bf s})=0$ then,
    there exists an atlas of Darboux coordinates
   on $E$ such that $\dv=D(x,\theta)$ in any coordinates
   from this atlas}.

  4) {\it In the case if all cohomological classes of differential forms
 of degree less than $n$ corresponding to the semidensity $\sqrt {\dv}$
  on even Lagrangian surface $L$
  are equal to zero also, then
   there exists an atlas of Darboux coordinates
   on $E^{n.n}$  adjusted to cotangent bundle structure of $L$ such
    that $\dv=D(x,\theta)$ in any coordinates
   from this atlas, i.e. $E$ can be identified with $\Pi T^*L$
   with volume form on $\Pi T^*L$ induced by volume form on $L$.

   This statement holds for another $(n.0)$-dimensional
    Lagrangian surface $\widetilde L$ if cohomological
    class of odd valued one-form $[\Psi](L,\widetilde L)$
     is equal to zero.}

   \smallskip

    This Corollary removes uncorrectness of the
     considerations about equivalence of conditions
     (\statementa), (\statementb) and (\statementc) in the Statement 1
     of Introduction, which was done in [\JMPd] and [\SchCMPa].
    On the other hand  some statements of this Proposition
   in non explicit way were contained in
     the statements of Lemma 4 and Theorem 5 of the paper
      [\SchCMPa].

                 \bigskip

  \centerline{\bf 5. Invariant densities on  surfaces}

  \medskip

 First we recall shortly the problem of construction of invariant densities
  in sympelctic (super)manifolds.
  Then we
   consider explicit formulae for
    the odd invariant semidensity on
    non-degenerate surfaces of codimension $(1.1)$
    embedded in an odd symplectic supermanifold $E$ provided with a
    volume form $\dv$ ([\CMP,\But]).
      We consider this semidensity
      as a kind of pull-back of
      semidensity  ${\bf s}$
      from the ambient odd symplectic supermanifold
       on embedded $(1.1)$-codimensional surfaces
        in the case if ${\bf s}=\sqrt{\dv}$.
       Using this construction for the semidensity
     $\Ds\sqrt\dv$ we will construct the new densities
     on embedded non-degenerated surfaces.

\smallskip

   In the case if we consider a volume form
   not only on the space (superspace) but
   on arbitrary embedded surfaces we come to
 the concept of densities on embedded surfaces.

  \noindent The density of weight $\s$ and rank $k$
  on embedded surfaces is a function
  $A(z,{\p z\over\p\zeta},\dots,{\p^k z\over\p\zeta\dots\p\zeta})$
   that is defined on parameterized surfaces $z(\zeta)$,
   depends on first $k$ derivatives of $z(\zeta)$
   and is multiplied on the $\s$-th power of determinant
(Berezinian) of surface reparametrization.

  A density of a weight $\s$ defines on every given
 surface $\s$-th power of volume form.
  Such a concept of density is very useful
  in supermathematics where the notion of differential forms
  as integration objects is ill-defined.
  It was elaborated by A.S.Schwarz, particularly
  for analyzing supergravity Lagrangians [\SchNucl, \Gayduk, \SchCMPb].

   In usual mathematics,
  for every $2k$-dimensional surface $C^{2k}$ embedded
   in a symplectic space,
    so called Poincar\'e-Cartan integral invariants
    (invariant volume forms on embedded surfaces)
     are given by the formula
     \def\cartan {5.1}
                       $$
   \int_{C^{2k}} \underbrace {\omega\wedge\dots\wedge \omega}_{k-{\rm times}}=
              \int
                 \sqrt
                   {
          \det
          \left
             (
        {\p x^\mu(\xi)\over\p\xi^i}
               \omega_{\mu\nu}
    {\p x^\nu(\xi)\over\p\xi^j}
              \right
            )
            }d^{2k}\xi\,,
                          \eqno (\cartan) 
          $$
 where a non-degenerate closed two-form
 $w=w_{\mu\nu}dx^\mu\wedge dx^\nu$ defines symplectic structure,
  and functions
 $x^\mu=x^\mu(\xi^i)$ define some parameterizations
 of the surface $C^{2k}$.

   In supermathematics
    one can consider even and odd symplectic structures
   on supermanifold generated by even and odd non-degenerate
   closed two-forms respectively [\Berezin,\leitb,\leites].

 In the case of an even symplectic supermanifold,
  the l.h.s. of (\cartan) is ill-defined but  the r.h.s.  of
  this formula  can be straightforwardly generalized,
   by changing determinant on the Berezinian
    (superdeterminant).
    The properties
    of the integral invariant do not change drastically.
   In particular one can prove  that
     the integrand in (\cartan)
 (the density of the weight $\s=1$ and of the rank $k=1$)
     is locally total derivative
       and all invariant densities
   on surfaces are exhausted by (\cartan)
   as well as in the case of usual symplectic structure
   [\Poin, \BarScha, \BarSchb].

    The situation is less trivial in the case of an odd  symplectic
     supermanifold.
   Formula (\cartan) cannot be generalized in this case
  because transformations preserving
  odd symplectic structure do not preserve any volume form . One can
 consider invariant densities only in an odd symplectic
  supermanifold provided with a volume form.

  The problem of the existence of invariant densities on non-degenerate
 surfaces embedded in an odd symplectic supermanifold
 provided with a volume form was studied in [\CMP, \But].
     In particularly it was proved that there are no invariant densities
   of the rank $k=1$ (except of the volume form itself),
    and invariant semidensity of the rank $k=2$
     that is defined on non-degenerate surfaces of the codimension
      $(1.1)$ was obtained. We briefly expose here its
      construction.

     The surface embedded in symplectic supermanifold
      is called  non-degenerate if the
     sympelctic structure of the supermanifold
      generates
     the symplectic structure on the embedded surface also,
     i.e. if the pull-back of the symplectic $2$-form
     on the surface is non-degenerate
      $2$-form. This symplectic structure
      on an
      embedded surface is called induced symplectic structure.

          Let $\{z^A\}$ be Darboux coordinates
    on an odd symplectic supermanifold
    $E=E^{n.n}$ provided
    with volume form $\dv=\rho(z)Dz$.
    It is convenient in this section to use for Darboux coordinates
    notations
    $z^A=(x^\mu,\theta_\mu),$ $(\mu=(0,i)=$ $(0,1,\dots,n-1)$,
    $i=(1,.\dots,n-1))$.
       Let $z(\zeta)$ be an arbitrary parametrization
       of an arbitrary non-degenerate
       surface of codimension $(1.1)$, embedded
        in $E$.
    ($\zeta=(\xi^i,\eta_j)$, $\xi^i$ and $\eta_j$
    are even and odd parameters respectively,
    ($i,j=1,\dots,n-1$)).
       The invariant semidensity of the rank $k=2$
      (depending on first and second derivatives of $z(\zeta)$)
       is given by the
      following formula [\CMP]:
  \def\odin {5.2}
                    $$
              A\left(z(\zeta),{\p z\over\p\zeta},
    {\p^2 z\over\p\zeta\p\zeta}\right)\sqrt{D\zeta}=
                    $$
            $$
               \left(
    \Psi^A{\p \log\rho(z)\over\p z^A}\,
  -\,\Psi^A\O_{AB}{\p^2 z^B\over\p\zeta^\a\p\zeta^\beta}
            \O^{\a\beta}(z(\zeta))
                  (-1)^
       {p(z^B(\zeta^\a+\zeta^\beta)+\zeta^\a}
               \right)
             \sqrt{D\zeta}\,,
                                    \eqno (\odin) 
                   $$
  where $\O_{AB}dz^A dz^B$ is the two-form defining the odd
  sympelctic
   structure on $E^{n.n}$
   and $\O^{\a\beta}$ is the tensor inversed to
   the two-form that defines induced symplectic structure on the
   surface.
   The vector field ${\bf \Psi}=\Psi^A{\p\over \p z^A}$
    is defined as follows:
   one have to consider
   the pair of vectors $(\bf H,\bf\Psi)$,
   $\bf H$ even and $\bf\Psi$ odd
   that are symplectoorthogonal to the surface
    and obey to the following conditions:
    \def\dva{5.3}
                    $$
      \O\left({\bf H,\bf\Psi}\right)=1,\,\,
        \O\left({\bf\Psi,\Psi}\right)=0
     \quad ({\rm symplectoorthonormality\,
     conditions})\,,
                     \eqno (\dva) 
           $$
\def\tri{5.4}
           $$
       d{\bf v}\left(\left\{{\p z\over\p\zeta}\right\},
       {\bf H},{\bf\Psi}\right)=1\quad
        \hbox{(volume form normalization conditions)}\,.
                           \eqno (\tri) 
                $$
These conditions fix uniquely the vector field ${\bf \Psi}$. (See
for  details [\CMP]).

 The explicit expression for this semidensity
  was calculated in [\But] in  terms
    of dual densities: If $(1.1)$-codimensional
   surface $C$ is given not by parametrization, but
    by the equations
      $f=0,\varphi=0$,
   where $f$ is an even function and $\varphi$ is an odd function
 then to the semidensity (\odin) there corresponds the
 dual semidensity:
 \def\four{5.5}
                   $$
                 {\tilde A}\Big\vert_{f=\varphi=0}=
           {1\over\sqrt{\{f,\varphi\}}}
                  \left(
                  \D f-
          {\{f,f\}\over 2\{f,\varphi\}}
                \D\varphi-
              {\{f,\{f,\varphi\}\}\over\{f,\varphi\}}-
                  {\{f,f\}\over 2\{f,\varphi\}^2}
              \{\varphi,\{f,\varphi\}\}
                       \right)\,.
                                       \eqno(\four) 
                        $$
  One can check that r.h.s. of (\four)
  restricted by conditions $f=\varphi=0$
 is multiplied by the square root of the corresponding
  Berezinian (superdeterminant) under the
   transformation
     $f\rightarrow$ $a f+\alpha\varphi$,
     $\varphi\rightarrow$ $\beta f+b\varphi$,
   which does not change the surface $C$ [\But].

  This invariant semidensity takes odd values.
 It is an exotic analogue of Poincar\'e--Cartan
 invariant: the corresponding density
 (the square of this odd semidensity)
 is equal to zero, so it cannot be  integrated
 nontrivially over surfaces.
      On the other hand this semidensity can be considered as an analog
      of the mean curvature of hypersurfaces in the Riemanian
    space [\CMP].

      This odd semidensity in an
      odd symplectic supermanifold
      is unique (up to multiplication
       by a constant) in the class of densities of the rank $k=2$
       that are defined on non-degenerated
       surfaces of arbitrary dimension
        [\But].
     This means that one have to search
      non-trivial integral invariants
      (invariant densities  of weight $\s=1$)
     in higher order derivatives (rank $k\geq 3$).
     Tedious calculations, which lead to
     the construction of the odd invariant semidensity
     in the papers [\CMP,\But]
     did not give hope to
    go further for finding them,
     using the technique used in these papers.

   Now we develop another approach
 rewriting the semidensity (\odin)
 straightforwardly via the semidensity $\sqrt\dv$
 on the ambient odd symplectic supermanifold $E=E^{n.n}$.

  Consider
  for every given non-degenerate surface $C$ of codimension
 $(1.1)$ embedded in odd symplectic supermanifold $E$
 Darboux coordinates
 such that in these Darboux coordinates
 the surface $C$ locally is given by equations
 \def\five{5.6}
                       $$
              x^0=\theta_0=0\,.
               \eqno (\five) 
              $$
We call these Darboux coordinates adjusted to the surface $C$.
 (The existence of Darboux coordinates obeying to
 these conditions can be proved
 using technique considered in Appendices 2 and 3).

  If $\{x^\mu,\theta_\nu\}$
are Darboux coordinates in $E$
 adjusted to the surface $C$, then
 $\{x^i,\theta_j\}$ are Darboux coordinates
 on the surface $C^{n-1.n-1}$  w.r.t. the induced symplectic structure
 ($\mu,\nu=0,\dots,n-1$, $i,j=1,\dots,n-1$).

Consider a semidensity (\odin) on arbitrary non-degenerated
surface $C=C^{n-1.n-1}$ of codimension $(1.1)$
 in Darboux coordinates (\five) adjusted to this
surface. Conditions of symplectoorthonormality in (\dva) give
that  ${\bf H}={1\over a}{\p\over \p
x^0}+\beta{\p\over\p\theta_0}$ and ${\bf\Psi}=a{\p\over\p
\theta_0}$, where $a$ is even and $\beta$ is odd. The condition
(\tri) of the volume form
 normalization gives that
                    $$
 a=\sqrt\rho\,{\rm Ber}^{1/2}
 \left(
 {\p(x^i,\theta_j)\over\p(\xi^i,\eta_j)}
 \right)\,,
             $$
 where a volume form $\dv$ is equal to $\rho(x,\theta) D(x,\theta)$
 and $\zeta=(\xi^i,\eta_j)$ are parameters
 ($x^0=\theta_0=0$, $x^i=x^i(\xi,\eta)$,
 $\theta_j=\theta_j(\xi,\eta)$).

 Hence the semidensity (\odin) on
a surface (\five) is reduced to
               $$
             A\left(z(\zeta),{\p z\over\p\zeta},
    {\p^2 z\over\p\zeta\p\zeta}\right)\sqrt{D\zeta}=
             a{\p \log \rho\over\p\theta_0}\sqrt{D\zeta}=
             2{\p \sqrt \rho\over\p\theta_0}\sqrt{D(x^i,\theta_j)}\,.
                    $$

   We come to the following statement

 {\bf Theorem}
   {\it To every semidensity ${\bf s}$ in the
   odd symplectic supermanifold $E$
   corresponds semidensity $\K({\bf s})$
   of an opposite parity defined on non-degenerated
   $(1.1)$-co\-di\-men\-si\-o\-nal
   surfaces embedded in this supermanifold.

   If semidensity ${\bf s}$ is given by expression
   ${\bf s}=s(x,\theta)\sqrt{D(x,\theta)}$
   in Darboux coordinates
    $\{x^\mu,\theta_\nu\}=$ $\{x^0,x^i,\theta_0,\theta_j\}$
   adjusted
    to given non-degenerate surface $C$
    of codimension $(1.1)$
    ($x^0\vert_C=$ $\theta_0\vert_C=0$),
    then semidensity $\K({\bf s})$
   on this surface
   in these Darboux coordinates
    is given by the following expression}:
    \def\theorem {5.7}
                          $$
   \A({\bf s})\big\vert_C={\p s(x^\mu,\theta_\nu)\over\p\theta_0}
                  \Big\vert_{x^0=\theta_0=0}
                  \sqrt{D(x^i,\theta_j)}\,,
                         \eqno (\theorem) 
                 $$

 where $D(x^\mu,\theta_\nu)$ is coordinate
 volume form on the supermanifold $E$
 and  $D(x^i,\theta_j)$ is coordinate
 volume form on the surface $C$.

 \smallskip

   The considerations above lead to the statement of this Theorem
 for semidensities related with a volume form
 on an odd symplectic supermanifold
   (${\bf s}=\sqrt{\dv}$),
  i.e. for even non-degenerate even semidensities.
  Continuity considerations lead to the fact that
  the formula (\theorem) is well-defined for an arbitrary
  semidensity
   e.g. for an odd semidensity,
 when the corresponding volume form is equal to zero.

 Alternatively one can prove this Theorem
checking in a same way as for (\deltadensdef)
 that the semidensity in r.h.s. of (\theorem)
is well defined.
 For example consider canonical transformation
  that has the following appearance in Darboux coordinates
  adjusted to surface $C$:
                  $$
                   \cases
                      {
     \t x^0=\t x^0(x^0,\theta_0)\,,\,
     \t\theta_0=\theta_0(x^0,\theta_0)\cr
   \t x^i=x^i\,,
     \t\theta_i=\theta_i\cr
           }
           $$
 One can see that these canonical transformations
  are exhausted by transformations
  $\t x^0=f(x^0)$, $\t\theta_0=\beta(x^0)+\theta_0/f_x$,
  $\t x^i=x^i,\theta_i=\theta_i$,
where $f(x)$ and $\beta(x)$ are even-valued and odd valued
functions on $x$ respectively. Hence for transformation of
adjusted coordinates $\t x^0=f(x^0)$, $\t\theta_0=\theta/f_x$.
Obviously r.h.s. of formula (\theorem)
 transforms as semidensity under this transformation.
 This is the central point of the construction (\theorem)
 and also of (\odin) (see for details [\CMP]).)

   We can consider a semidensity $\K({\bf s})$ in (\theorem)
    as a kind of pull-back
   of semidensity ${\bf s}$ on $C$, but this construction does not
   obey to condition of transitivity for pull-back:
   consider arbitrary $(k.k)$-dimensional non-degenerate surface
   embedded in $E^{n.n}$ and include it in
   a flag of non-degenerated surfaces:
   \def\flag{5.8}
                            $$
       Y^{k.k}\hookrightarrow Y^{k+1.k+1}
       \dots\hookrightarrow Y^{n-1.n-1}\hookrightarrow E^{n.n}\,.
                 \eqno (\flag)
                            $$
    then
    one can consider semidensity $\A(\dots \A({\bf s})\dots)$
    on $Y^{k.k}$ corresponding to the semidensity
    ${\bf s}$  depending on flag (\flag).

The statement of Theorem allows us to construct semidensity on
embedded surfaces via odd semidensities on the ambient
supermanifold, which cannot be yielded from volume forms.

  In an odd sympelctic supermanifold provided with a volume form $\dv$
  on $(1.1)$-co\-di\-men\-si\-onal non-degenerate surfaces
  except an odd semidensity $\K(\sqrt\dv)$, that is nothing
    but semidensity (5.2),
  one can consider also an even
 semidensity $\K(\Ds\sqrt\dv)$ corresponding
 to an odd semidensity $\Ds\sqrt {\dv}$.
  The semidensity $\K(\Ds\sqrt\dv)$ cannot be represented
 (5.2)-like because the square of the odd semidensity
 $\Ds\sqrt\dv$ is equal to zero.

  We note that for semidensities $\K(\Ds{\bf s})$
  and $\K({\bf s})$ for arbitrary $(1.1)$-codimensional surface
  $C$
  the following condition is obeyed:
  \def\posleflaga {5.9}
                  $$
    \K(\Ds{\bf s})\Big\vert_C=
    -{\widetilde {\Ds}}\K({\bf s})\Big\vert_C\,,
                    \eqno (\posleflaga)
                   $$
  where   ${\widetilde {\Ds}}$ is $\Ds$-operator on
  surface $C$ w.r.t. induced symplectic structure.
  This relation can be immediately checked in Darboux coordinates
  (\five) adjusted to the surface $C$.

The semidensities $\K(\sqrt\dv)$ and $\K(\Ds\sqrt\dv)$
 can be integrated
over Lagrangian
 subsurfaces in $C$, according (\volumemap).

 On the other hand one can  consider the following
 non-trivial densities of weight $\s=1$ constructed
 via the semidensities  $\A(\sqrt\dv)$ and
 $\A(\Ds\sqrt\dv)$:
\def\newdensities {5.10}
                  $$
          P_0=\A^2(\Ds\sqrt\dv)\,\quad{\rm and}\quad
          P_1=\A(\sqrt\dv)\A(\Ds\sqrt\dv)\,.
               \eqno (\newdensities) 
             $$
The density  $P_0$ takes even values, the density $P_1$ takes odd
values. In general case these densities give non-trivial
integration objects (volume forms) over non-degenerated
$(1.1)$-codimensional surfaces embedded in an odd symplectic
supermanifold with volume form $\dv$.

The densities $P_0$ and $P_1$ have rank $k=4$
 (i.e. depend
on derivatives of the parametrization $z(\zeta)$ up fourth order).
It follows from the fact that according to (\posleflaga) the
semidensity $\K(\Ds\sqrt\dv)$ has the rank $k=4$,
 because the semidensity $\K(\sqrt\dv)$
 has the rank $k=2$.
 This is hidden in representation (\theorem), where the function
  $\rho(z)$ corresponding to the volume
 form in adjusted coordinates depends
 non-explicitly on
 derivatives of surface parametrization $z(\zeta)$.

Finally we consider a simple example of these constructions and
their relations with differential forms.

  Let $E^{3.3}$ be a superspace associated to
 cotangent bundle of $3$-dimensional space $E^3$, $E^{3.3}=\Pi T^*E^3$.
 We assume that coordinates ${x^0,x^1,x^2}$ are globally defined on $E^3$.
 We consider on $E^3$ the differential form
    $$
    w=-dx^0\wedge dx^1\wedge dx^2+b_0 dx^0+b_1 dx^1+b_2 dx^2\,.
    $$
    According to (3.4)
    a semidensity
          $$
  {\bf s}=\tau^{\#}(w)=
  (1+b_0\theta_1\theta_2+b_1\theta_2\theta_0+
  b_2\theta_0\theta_1)
  \sqrt{D (x^0,x^1,x^2,\theta_0,\theta_1,\theta_2)}
            $$
 in $\Pi T^*E$ corresponds to this differential form.
Let $C$ be a surface in $E^{3.3}$ that is defined by equations
 $x^0=\theta_0=0$ and $E^2$ is subspace in $E^3$ defined by equation $x^0=0$.
  $(2.2)$-dimensional superspace $C$ provided with coordinates
  ${x^1,x^2,\theta_1,\theta_2}$ can be identified with
  superspace $\Pi T^*E^2$ associated with cotangent bundle $T^*E^2$
  of subspace $E^2$.

 Then the  value of the odd semidensity $\K({\bf s})$ on
 $C=\Pi T^*E^2$ is equal
 to $(b_2\theta_1-b_1\theta_2)
 \sqrt{D (x^1,x^2,\theta_1,\theta_2)}$.
 This semidensity corresponds to differential form $b_1dx^1+b_2dx^2$,
 the pull-back of $w$ on $E^2$.
 The value of even semidensity $\A(\Ds{\bf s})$ on $C$ is equal
 to $(\p_2b_1-\p_1 b_2)\sqrt{D (x^1,x^2,\theta_1,\theta_2)}$.
 This corresponds to differential form $d(b_1dx^1+b_2dx^2)$
 $=(\p_2b_1-\p_1 b_2)dx^1\wedge dx^2$, the pull-back
 of $dw$ on $E^2$.

 The even density (volume form) on $M$ is equal to
  $P_0=(\p_2b_1-\p_1 b_2)^2D (x^1,x^2,\theta_1,\theta_2)$
 and the odd density $P_1$
 $=(\p_2b_1-\p_1 b_2)(b_1\theta_2-b_2\theta_1)
   D (x^1,x^2,\theta_1,\theta_2)$.

\bigskip

\centerline {\bf 6.Discussion}

\medskip

 The definition (\deltadefgen) of $\D$-operator is  well-defined
 not only for symplectic supermanifold but
 for Poisson supermanifold also (provided with a volume form)
 even if corresponding Poisson bracket is degenerate.
 Is it possible to define $\Ds$-operator on semidensities
 in odd Poisson supermanifold?

 We note also that from relations (\deltaprop)
 it follows that one can
 express odd Poisson bracket via operator $\D$. Moreover
 every second
 order odd differential operator $\hat A$ on functions on
 supermanifold without the derivateveless
 term and obeying to condition $\hat A^2=0$ defines
 Poisson structure via relations (\deltaprop).
  (This approach was elaborated by I.A. Batalin and I.V. Tyutin [\tyutina].)
  It is interesting to consider the analogue of this construction
 for differential operators on semidensities.

In an odd symplectic supermanifold integration theory on surfaces
interplays with symplectic geometry.
 Using our approach one can consider integrands
 (differential forms) of
  functionals in terms of semidensities
  corresponding to these
 differential forms.
 In this case
  symmetry transformations of these functionals
  are not exhausted only by transformations induced by
  diffeomorphisms  of underlying space. General canonical
  transformations of supermanifold induce
  mixing of corresponding differential
  forms with different degrees.

In Sections 3 and 4 we investigated relations between
semidensities on an odd symplectic supermanifold and
 differential forms on even Lagrangian
 surfaces.
    It is interesting
   to generalize these results
  on the case when Lagrangian surface is
  $(n-k.k)$-dimensional for $k\not=0$
  on the base of analysis of arbitrary Lagrangian surfaces
  performed in the paper [\SchCMPa].
   In this case there exist analogies of differential
   forms considered as integration object.
   (See [\Bernsteina, \Bernsteinb], [\BarScha, \BarSchb] and
   for a more detailed analysis [\Vora, \Vor].)
   For example if $\Pi T^*M$
  is $(r+s.r+s)$-dimensional supermanifold associated with cotangent
  bundle of $(r.s)$-dimensional supermanifold $M$
  then considering analogue of formula (3.4a) we
  obtain relations between
  semidensities in $\Pi T^*M$ and so called
  pseudodifferential forms: functions
  on supermanifold $\Pi TM$. Pseudodifferential
  forms are well-defined integration objects
  over supermanifold and embedded surfaces
  [\Bernsteina, \Bernsteinb,
   \BarScha, \BarSchb, \Vora, \Vor].
  In this way it is possible to come to
  an analogue of the map (\volumemap)
 (see [\SchCMPa], [\JMPd]). It is interesting  to
 construct analogues of maps (\mapcoh) and (\cohomologymap)
 for $(n-k.k)$-dimensional Lagrangian
 surfaces.

  We note that our considerations
  in subsections 4.2 and 4.3 overlap partially  with
  some results of the paper [\SchCMPa].
  But
  we perform analysis in terms of semidensities
  where calculus analogous to the calculus of differential
  forms arises.  This makes our considerations more exact
  that considerations in terms of volume forms performed in [\SchCMPa].
   In particularly this leads to exact statements
  in Corollary 2. We note also that considering $\Lambda$-points
  formalism for supermanifold we come to the difference between
  supergroups $Can_0(E)$ and $Can_H(E)$.

 We hope that considerations presented
  in Section 5 of this paper
 can be generalized for constructing densities
 depending on higher order derivatives
  for surfaces of arbitrary  dimension
  embedded in an odd symplectic supermanifold provided
   with a volume form
 and for finding the complete set of local invariants of this
 geometry. In particular from considerations which lead to Theorem follows
that if $k(p)$ is rank of
 non-trivial invariant densities
  on non-degenerated surfaces of codimension $(p.p)$,
  then $k(2)\geq 5$ and $k(p+1)>k(p)$.

 In [\CMP] some relations
 of semidensity (\odin, \theorem)
 with mean curvature in Riemanian geometry
  were indicated. It is interesting to analyze
  these relations in terms of
 geometrical interpretations of semidensities presented in this paper.

 The densities presented in formula (\newdensities)
 are needed to be investigated more in details. Particularly
  one have to present explicit formulae for them and
  consider
  the corresponding functionals over surfaces.
  These functionals are equal to zero in the special case
   if the volume form in the ambient odd
   symplectic supermanifold obeys to BV-master equation.
 Do Euler-Lagrange motion equations
 for these functionals in a general case
 are equal identically to zero,
 as for usual Poincare-Cartan integral invariants (\cartan)?

\smallskip

 Results presented in
 Section 5 strongly indicate  that
  there exists non-trivial geometry
  in an odd symplectic supermanifold
  provided with a volume form.

\bigskip

              \centerline{\bf  Acknowledgment}
              \medskip
 This work is highly stimulated by
very illuminating discussions with S.P. Novikov during my talk on
his seminar in Moscow in August 1999. I am deeply grateful to him.

I want to express my deep gratitude to my teacher A.S.Schwarz, and
to I.A.Batalin,

\noindent V.M.Buchstaber, I.V.Tyutin for encouraging me to do
this work.

I am very grateful to T. Voronov.
 Continuous discussions with him
were very useful in the final part of this work.
 He also did enormous work by reading a draft of
 this paper and giving many
 valuable advises.

I am deeply grateful for hospitality and support of the Abdus
Salam International Centre for Theoretical Physics in Trieste and
Max-Plank-Institut f\"ur Mathematik in Bonn, where I began a work
on this paper.

  The work was partially supported by grant EPSRC GR/N00821.

\bigskip

\centerline {\bf  Appendix 1. $\L$-points on Supermanifolds}

\medskip

  Let $\left[\{x^i_{(\a)}\}\right]$ be a smooth atlas of
coordinates
 on $m$-dimensional manifold $M^m$,
 where coordinates $\{x^i_{(\a)}\}$ are defined on domain
 $U_\a$ and
 $x^i_{(\a)}=\Psi^i_{\a\beta}(x_{(\beta)})$ are transition functions.
  Consider an atlas
 $\left[\{x^i_{(\a)},\theta^j_{(\a)}\}\right]$, where
 odd variables $\{\theta^j_{(\a)}\}$ $(j=1,\dots,n)$ are
 generators of Grassmann algebra
 and transition functions
                      $$
    \cases
      {
 x^i_{(\a)}=
 \widetilde\Psi^i_{\a\beta}(x_{(\beta)},\theta_{(\beta)})\cr
\theta^j_{(\a)}=
 \Phi^j_{\a\beta}(x_{(\beta)},\theta_{(\beta)})\cr
               }
               \eqno ({\rm Ap}1.1)
               $$
               obey to the following properties:

\noindent 1) they are parity preserving, i.e.
   $p(\widetilde\Psi_{\a\beta})=0, p(\Phi_{\a\beta})=1$, where
   $p(x^i)=0, p(\theta^j)=1$,

 \noindent 2)
 $\widetilde\Psi_{\a\beta}(x_{(\beta)},\theta_{(\beta)})\big\vert_{\theta^j=0}=
  \Psi_{\a\beta}(x_{(\beta)})$
 and ${\p\Phi^j/\p\theta^{i}_{(\beta)}}$ are inverting matrices.

 Coordinates $\{x^i_{(\a)},\theta^j_{(\a)}\}$ define $(m.n)$-dimensional
 superdomain
 $\hat U^{m.n}_{(\a)}$ with underlying domain $U^m_{(\a)}$.
 Pasting formulae (Ap1.1) define $(m.n)$-dimensional
 supermanifold with underlying manifold
 $M^m$. In this definition of \superspace which belongs to F.Berezin
 and D.Leites (see [\Berezin] and [\leites]) a supermanifold "has
 no points".

 If $E$ is supermanifold and $\L$ is an arbitrary Grassmann
 algebra one can construct a set $E_\L$ of $\Lambda$-points
 of supermanifold $E$. For example if $E^{m.n}$ is superdomain
 with underlying domain $M^m$,
 we define $E_\L$ as a set of rows $(a^1,\dots,a^m,\a^1,\dots,\a^n)$,
 where $a^1,\dots,a^m$ are arbitrary even elements and
 $\a^1,\dots,\a^n$ are arbitrary odd elements of Grassmann algebra $\L$
 and $(m(a^1),\dots,m(a^m))\in M^m$, where $m$ is a standard homomorphism
 of $\L$ on $I\!R$. A map of superdomains generates a map of corresponding
 sets of $\L$-points. Thus one comes to definition of a set $E_\Lambda$
 for arbitrary supermanifold $E$.
   To every parity preserving homomorphism
  $\rho\colon$ $\L\rightarrow\L^\prime$ of Grassmann algebras one
  can naturally assign a map $\tilde\rho_{_E}\colon$
  $E_\L\rightarrow E_{\L^\prime}$. If
  $\rho\colon\quad $ $\L\rightarrow\L^\prime$  and
  $\rho^\prime\colon$ $\L^\prime\rightarrow\L^{\prime\prime}$
  are two parity preserving homomorphisms, then
  $\widetilde {(\rho\circ\rho^\prime)_{_E}}=
  \tilde{\rho_{_E}}\circ\tilde{\rho_{_E}^\prime}$.
  Supermanifold can be considered as functor on the category of
  Grassmann algebras taking values in category of sets.

  This definition of supermanifolds
   is used in the paper.
   It was suggested and widely used by A.S. Schwarz [\SchCMPb].
   It  makes possible to use
  a language of "points" and is more convenient for supergeometry
  and in applications in theoretical physics
  \footnote{$^*$}
  {A possibility to use a language of $\L$-points
  was noted by D.Leites in [\leites].}.

 In terms of $\L$-points one can easy to generalize the
 standard geometrical definitions on supercase [\SchCMPb].
 For example

 1. A map $F$ from
\superspace $E$ in \superspace $N$ can be considered as a functor
from category
   $\{\L\}$ of Grassmann algebras to category
   $\{F_\L\}$, where for every Grassmann algebra $\L$
  $F_\L$
  is a map from the set
 $E_\L$ to the set $N_\L$
  such that
 $F_{\L^\prime}\circ\tilde{\rho_{_E}}=
 \tilde{\rho_{_N}}\circ F_{\L}$ for every parity preserving
 homomorphism $\rho\colon$ $\L\rightarrow\L^\prime$.

 2. The action of supergroup $G$ on supermanifold $E$
 can be considered as a functor that assigns to every Grassmann algebra the
 pair $[G_\Lambda, E_\Lambda]$ where $G_\Lambda$ is
 a group of $\L$-points of supergroup $G$,
 that acts on the set $E_\Lambda$
 of $\Lambda$-points of supermanifold $E$.

\bigskip

   \centerline {\bf Appendix 2. A simple proof of Darboux Theorem}
\centerline {\bf for odd symplectic structure}
\medskip

   Using nilpotency of odd variables
  one can directly prove Darboux theorem for an odd symplectic
  supermanifold presenting finite recurrent procedure
  for constructing
  Darboux coordinates starting from arbitrary coordinates.

 Let $\{\quad,\quad\}$ be odd non-degenerated
    Poisson bracket (\bracket)
     corresponding to the symplectic structure.
    According to (2.1) for arbitrary two functions $f$ and $g$
                  $$
                  \matrix
                    {
     \{f,g\}={\p f\over\p x^i}\{x^i,x^j\}{\p g\over\p x^j}+
     {\p f\over\p x^i}\{x^i,\theta_j\}{\p g\over\p \theta_j}+
    (-1)^{p(f)+1}
    {\p f\over\p\theta_i}\{\theta_i,x^j\}{\p g\over\p x^j}\cr
                      +
    (-1)^{p(f)+1}
    {\p f\over\p\theta_i}\{\theta_i,\theta_j\}{\p g\over\p \theta_j}\cr
                   }
               \eqno ({\rm Ap2}.1)
                   $$
        and Jacoby identities (\jacoby) are obeyed.

 For given arbitrary coordinates $\{x^1,\dots,x^n,\theta_1,\dots,\theta_n\}$
  denote by
                 $$
        E^{ij}(x,\theta)=\{x^i,x^j\}\,,
       F_{ij}(x,\theta)=\{\theta_i,\theta_j\}\,,
A^i_j(x,\theta)=\delta^i_j+P^i_j(x,\theta)=\{x^i,\theta_j\}\,.
                                  \eqno ({\rm Ap}2.2)
                     $$
From definition of symplectic structure it follows that
 $E^{ij}=E^{ji}$, $F_{ij}=-F_{ji}$ are odd-valued matrices
 taking values in Grassmann algebra $\L$
 and $A^i_j(x,\theta)$ is even non-degenerate
 matrix taking values in Grassmann algebra $\L$.
In Darboux coordinates matrices $E^{ij}$, $F_{ij}$ and $P^i_j$
have to be equal to zero.

First of all we note that in the case if for coordinates
$\{x^i,\theta_j\}$ the conditions
                         $$
            E^{ik}(x,\theta)=0\,,\quad P^i_k(x,\theta)=0
                              \eqno ({\rm Ap}2.3)
                          $$
  are obeyed then
  Jacoby identities
  $\{x^m,\{\theta_i,\theta_j\} \}+\hbox{cycl. permut.} =0$
  imply that $F_{ij}$ do not depend on $\theta$
  and Jacoby identities
  $\{\theta_i,\{\theta_j,\theta_m\}\}+
  \hbox{cycl. permut.}=0$
  imply the condition $\p_i F_{jm}(x)+\p_j F_{mi}(x)+\p_m F_{ij}(x)=0$.
  (In other words two-form $F_{ij}(x)dx^i\wedge dx^j$ is closed).
  Locally it means that there exist functions $A_i(x)$
  such that $F_{ij}(x)=\p_i A_j(x)-\p_j A_i(x)$.
   Under transformation $\theta_i\rightarrow \theta_i+A_i(x)$,
  $F_{ij}(x)$ transform to zero also and we come to Darboux coordinates.

\def\M{{\cal M}}
  Thus we have to find transformation from arbitrary coordinates
   to new coordinates such that
  in new coordinates conditions (Ap2.3) will be obeyed.

   Consider a set $\M$ of all coordinates $\{x^i,\theta_j\}$.
   Denote by $\M_{(p.q)}$ a subset of $\M$ such that for coordinates
   $\{x^i,\theta_j\}$ belonging to the subset $\M_{(p.q)}$
    the following conditions are
   obeyed for matrices $E^{ij}(x,\theta)$ and $P^i_j(x,\theta)$ in (Ap2.2):
                            $$
   E^{ij}(x,\theta)=O(\theta^p)\,,\quad
    P^i_j(x,\theta)=O(\theta^q)\,.
                        \eqno ({\rm Ap}2.4)
                   $$
$\M_{0.0}=\M$ and condition $\{x^i,\theta_j\}\in\M_{n+1.n+1}$
means that relations (Ap2.3) are obeyed for these coordinates,
because $\theta_{i_1}\dots\theta_{i_k}=0$ if $k\geq n+1$.

 Consider four maps $\F_1,\F_2,\F_3,\F_4$ defined
  on the set $\M$ of coordinates, such
  that these maps obey to the following conditions:
                           $$
                           \matrix
                             {
  \hbox{$\F_1$ maps $\M_{r.0}$ in $\M_{r.1}$  for $r=0,1,\dots$}\quad,\cr
   \hbox{$\F_2$ maps $\M_{0.1}$ in $\M_{1.0}$ }\,,\cr
  \hbox{$\F_3$ maps $\M_{r.1}$ in $\M_{r+1.1}$ for $r\geq 1$}\,,\cr
  \hbox{$\F_4$ maps $\M_{n+1.r}$ in $\M_{n+1.r+1}$
              for $r\geq 1$}\,.\cr
                           }
                              \eqno ({\rm Ap}2.5)
                           $$
Provided conditions (Ap2.5) are obeyed
   the map $\F_4^n\circ\F_3^n\circ \F_1\circ\F_2\circ\F_1$
transforms arbitrary coordinates to coordinates that belong to
subset $\M_{n+1.n+1}$, i.e. conditions (Ap2.3) are obeyed for
transformed coordinates.

 Now we present maps $\F_1,\F_2,\F_3,\F_4$
 obeying conditions (Ap2.5).

1. Definition of the  map $\F_1$:
                       $$
     \F_1(\{x^i,\theta_j\})=\{\t x^i,\t\theta_j\},\quad{\rm where}\quad
             \t x^i=x^i, \t\theta_j=\theta_m(A^{-1})_j^m\,,
                    \eqno  ({\rm Ap}2.6)
                    $$
where matrix $A^{-1}$ is inverse to the matrix $A$ defined by
relations (Ap2.2) for coordinates $\{x^i,\theta_j\}$. It is easy
to see from (Ap2.1) that map (Ap2.6) obeys condition (Ap2.5).

2. Definition of the map $\F_2$:
                    $$
                    \F_2(\{x^i,\theta_j\})=\{\t x^i,\t\theta_j\},
                    \quad{\rm where}\quad
              \t x^i=x^i-\theta_m R^{mi},\quad
              \t\theta_j=\theta_j\,,
                    \eqno  ({\rm Ap}2.7)
                          $$
where symmetrical odd-valued matrix $R$ is solution to matrix
equation
                    $$
                   2R+ RFR=E\,,\quad (R^{ij}=R^{ji})
                   \eqno ({\rm Ap2}.8)
                      $$
and matrices $E$  and $F$ for coordinates $\{x^i,\theta_j\}$ are
defined by (Ap2.2).

 The solution to this equation is well-defined because elements of
symmetric matrix $E$ and antisymmetric matrix $F$
 take odd values in Grassmann algebra $\L$.
$R$ is given by finite power series $R={E\over 2}-{EFE\over
8}+\dots$ containing less than $[{n^2\over 2}]$ terms.
 One can present explicit solution to equation (Ap2.8):
                 $$
       R=\sum_{k=0}^{{n(n-1)\over 2}} c_k(EF)^kE,\quad {\rm where}
               \quad\sum_{k=0}^{\infty}
                c_kt^k={\sqrt {1+t}-1\over t}\,.
                \eqno  ({\rm Ap}2.9)
                  $$

Now it follows from (Ap2.1) and (Ap2.8) that under transformation
(Ap2.7) matrix $E^{ij}=\{x^i,x^j\}$ transforms to the matrix $\t
E^{ij}=\{\t x^i,\t x^j\}$ such that
               $$
   \t E^{ij}=E^{ij}-2R^{ij}-R^{im}F_{mk}R^{mj}+O(\t\theta)=O(\t\theta)\,,
                               \eqno ({\rm Ap}2.10)
               $$
 if coordinates $\{x^i,\theta_j\}$ belong to $\M_{0.1}$
  (i.e. $A^i_j=\delta^i_j+O(\theta)$)
  and matrix $R$ obeys to equation (Ap2.8).
  Hence map (Ap2.7) obeys condition (Ap2.5).

3. Definition of the  map $\F_3$:
                    $$
                    \F_3(\{x^i,\theta_j\})=\{\t x^i,\t\theta_j\},
                    \quad{\rm where}\quad
             \t x^i=x^i-\theta_m
             \int_0^1\tau E^{mi}(x,\tau\theta)d\tau,\quad
             \t\theta_j=\theta_j\,.
                    \eqno  ({\rm Ap}2.11)
                          $$
  From (Ap2.1) it follows that
 transformation (Ap2.11) maps
 $\M_{r.1}$ in $\M_{r.1}$ if $r\geq 1$.
 Matrix $E^{ij}(x,\theta)$
  transforms to matrix
                 $$
                     E^{ij}-{2E^{ij}\over r+2}+
     {1\over r+2}\left(\theta_m{\p E^{mj}\over \p \theta_i}+
     (i\leftrightarrow j)\right)+O(\theta^{r+1})\,.
     \eqno ({\rm Ap}2.12)
                            $$
 On the other hand from Jacoby identity (\jacoby):
 $\{x^i\{x^j,x^m\}\}+\{x^j\{x^m,x^i\}\}+\{x^m\{x^i,x^j\}\}$
 $=0$
  and from (Ap2.1)
  it follows that
               $$
      \theta_m{\p E^{mj}\over \p \theta_i}+
     (i\leftrightarrow j)=
          -\theta_m {\p E^{ij}\over\p\theta_m}+O(\theta^{r+1})=
                 -rE^{ij}(x,\theta)+O(\theta^{r+1})\,.
                 \eqno ({\rm Ap}2.13)
                 $$
   Hence (Ap2.12) is equal to zero up to $O(\t\theta^{r+1})$
   and condition (Ap2.5) is obeyed for transformation
   (Ap2.11).

4. Definition of the map $\F_4$:
                    $$
                    \F_3(\{x^i,\theta_j\})=\{\t x^i,\t\theta_j\},
                    \quad{\rm where}\quad
             \t x^i=x^i,\,\,
             \t\theta_j=\theta_j-\theta_m
             \int_0^1 P^m_j(x,\tau\theta)d\tau\,.
                    \eqno  ({\rm Ap}2.14)
                    $$
 We prove that (Ap2.10) maps $\M_{n+1.r}$ in $\M_{n+1.r+1}$ analogously
 to the proof for (Ap2.11).
  Suppose that coordinates $\{x^i,\theta_j\}$ belong to
 $\M_{n+1.r}$ ($r\geq 1$). Then
 transformation

 \noindent (Ap2.14) maps matrix $P^i_j(x,\theta)$
 to matrix
                        $$
            P^i_j-{P^i_j\over r+1}+
            {\theta_m\over r+1}
         {\p P^m_j\over\p\theta_i}+O(\theta^{r+1})=
                   P^i_j-{P^i_j\over r+1}-
                    {\theta_m\over r+1}
         {\p P^i_j\over\p\theta_m}+O(\theta^{r+1})=
                O(\theta^{r+1})  \,,
                 $$
because of Jacoby identity $\{x^i,\{x^m,\theta_j\}\}+
\{x^m,\{x^i,\theta_j\}\}+ \{\theta_j,\{x^i,x^m\}\}=0$. Hence
condition (Ap2.5) is obeyed for
 transformation (Ap2.14).

\vfill\eject

\medskip
\centerline {\bf Appendix 3. Hamiltonians of adjusted
                         canonical transformations}

\medskip

 In this Appendix we prove that for any given
 adjusted canonical transformation  $\{x^i,\theta_j\}$
 $\rightarrow$
 $\{\t x^i,\t\theta_j\}$
 (\cantransform a) there exists time-independent Hamiltonian
 $Q(x,\theta)$ that generates this
 transformation via differential equations (\hamiltoniandeformation)
 and this Hamiltonian is defined uniquely by the condition
                        $$
 Q(x,\theta)=Q^{ik}\theta_i\theta_k+\dots\,,\quad
   {\rm i.e.}\,\,Q=O(\theta^2)\,.
                                       \eqno ({\rm Ap}3.1)
                      $$

\def\Ocal {{A}}

For every Hamiltonian (odd function) $Q(x,\theta)$ obeying
condition (Ap3.1) consider one-parametric family of functions
(Darboux coordinates)
 $\{y^i(t),\eta_j(t)\}$ ($i,j=1,\dots,n$)
 that are solution to differential equation (\hamiltoniandeformation):
                   $$
                    \cases
                     {
 {d y^i(t)\over dt}= \{Q(y,\eta),y^i\}=
 -{\p Q(y,\eta)\over \p\eta_i}\,,\cr
{d \eta_j(t)\over dt}= \{Q(y,\eta),\eta_j\}= {\p Q(y,\eta)\over\p
 y^i}\,,\cr
              }
              \quad
              (0\leq t\leq 1)\,,
                        \eqno ({\rm Ap}3.2)
                    $$
with initial conditions
                    $$
y^i(t)\big\vert_{t=0}=x^i,\,
  \eta_i(t)\big\vert_{t=0}=\theta_i\,.
                   $$

It is easy to see from explicit expression (\darbouxtheorem) for
odd Poisson bracket that if $\{ x^i,\theta_j\}$ and $\{\t
x^i,\t\theta_j\}$ are Darboux coordinates such that $\t x^i=x^i$
and $\t\theta_j=O(\theta)$ then $\t\theta_j=\theta_j$ also. Hence
every adjusted canonical transformation
$\{x^i,\theta_j\}\rightarrow$
 $\{\t x^i,\t\theta_j\}$ is uniquely defined
 by functions $\{f^i(x,\theta)\}$
  that obey to conditions:
                     $$
     \{x^i+f^i(x,\theta),x^j+f^j(x,\theta)\}=0\,\,{\rm and}\quad
               f^i(x,\theta)\in O(\theta)\,.
                          \eqno ({\rm Ap}3.3)
                     $$
Statement 3 of Lemma 1 follows from the Lemma:

{\bf Lemma 3}
   {\it For every set of functions $\{f^i(x,\theta)\}$ ($i=1,\dots,n$)
   obeying conditions} (Ap3.3) {\it there exists unique Hamiltonian
   $Q$ obeying condition} (Ap3.1)
    {\it such that functions $\{y^i(t)\}$ solutions to differential
   equation} (Ap3.2)
    {\it obey conditions $y^i(t)\vert_{t=1}=x^i+f^i(x,\theta)$
   ($i=1,\dots,n$).}

   Prove this Lemma.

 Consider a ring $A$ of functions on coordinates $(x^1,\dots,x^n$,
  $\theta_1,\dots,\theta_n)$.
  (As always functions take values in an arbitrary
  Grassmann algebra $\L$.
  Consider in $A$ the following gradation:
  $A_{(p)}$ is a space of functions that are linear combinations
   of $p$-th order monoms on variables $\{\theta_1,\dots,\theta_n\}$:
   $f\in A_{(p)}$ iff $\sum_k\theta_k{\p f\over\p\theta_k}=pf$.
   $A_{(p)}=0$ for $p\geq n+1$.
 For every function $f\in A$ we denote by $f_{(p)}$
 its component in $A_{(p)}$:
  $f=f_{(0)}+f_{(1)}+\dots+f_{(n)}$.
   It is evident that for canonical Poisson bracket
   (\darbouxtheorem)
                     $$
                     \{f,g\}_{(p)}=
         \sum_{i=0}^n \{f_{(i)},g_{(p+1-i)}\}\,.
                       \eqno ({\rm Ap}3.4)
                       $$
  Consider also a corresponding filtration:
              $$
               0=A^{(n+1)}\subset A^{(n)}\subset
 \dots\subset A^{(1)}\subset A^{(0)}=\Ocal\,,
              $$
 where $A^{(p)}=\oplus_{k\geq p} A_{(k)}$.


  We denote by $\Ocal^+$ ($\Ocal^-$) a subspace of even-valued
  (odd valued) functions in $\Ocal$.
   Respectively we denote by $\Ocal^{\pm}_{(k)}=\Ocal_{(k)}\cap\Ocal^{\pm}$
  and $\Ocal^{\pm(k)}=\Ocal^{(k)}\cap\Ocal^{\pm}$.

  We note first that condition (Ap3.1) implies that
  solutions to equations (Ap3.2) are well defined.
  Indeed consider arbitrary function $\varphi(x,\theta)$,
  odd Hamiltonian $Q\in A^{-(2)}$ and differential equation
  $\dot\varphi=\{Q,\varphi\}$. Projecting this differential
  equation on the subspace $A_{(p)}$ we come  using (Ap3.4) to
  equations $\dot\varphi_{(p)}=\{Q_{(p+1)},\varphi_{(0)}\}+\dots+
  \{Q_{(2)},\varphi_{(p-1)}\}$. Function $\varphi_{(0)}$
  does not depend on $t$ ($\dot \varphi_0=0$)
   and these equations can be solved
  recurrently:
                      $$
  \varphi_{(p)}\big\vert_{t=a}=\varphi_{(p)}\big\vert_{t=0}+
  a\{Q_{(p+1)},\varphi_{(0)}\}+\dots,
                \eqno ({\rm Ap}3.5)
                $$
where we denote by dots terms depending on $Q_{(2)},\dots,Q_{(p)}$
and functions $\varphi_{(0)},\varphi_{(1)}\vert_{t=0},\dots,$
$\varphi_{(p-1)}\vert_{t=0}$.

 \def\N{{\cal N}}
   \def\U {{\cal U}}
 Denote by $\N$ a space  of sets of even-valued functions
 $\{f^i(x,\theta)\}$
  ($i=1,\dots,n$) such that these functions obey to condition (Ap3.3).
  Consider a map that assigns to every Hamiltonian $Q\in A^{-(2)}$
  the solutions $\{y^i(t)\vert_{t=1}\}=x^i+f^i(x,\theta)$
   to differential equations
  (Ap3.2). Thus we define map $\U\colon$ $A^{-(2)}\rightarrow\N$.
   Relations (Ap3.5) for $\varphi=x^i$ imply that
                $$
      f^i_{(p)}=-{\p Q_{p+1}\over\p\theta_i}+
      \hbox{terms depending on
            $Q_{(2)},\dots,Q_{(p)}$}\,.
            \eqno ({\rm Ap}3.6)
              $$

   Consider also a following  map
   $\delta\colon$ $\N\rightarrow A^{-(2)}$
   such that  for every $\{f^i\}\in\N$
                       $$
   \delta(\{f^i(x)\})=
   -\sum_{i=1,p=1}^n \theta_i
                 {f^i_{(p)}(x,\theta)\over p+1}
                    =
    -\sum_{i=1}^n \theta_i\int_0^1 f^i(x,\tau\theta)d\tau\,.
                \eqno ({\rm Ap}3.7)
               $$
 From condition (Ap3.3) for functions $\{f^i\}$ and (\darbouxtheorem)
    it follows that
                    $$
   f^i=-{\p \t Q\over \p\theta^i}+
   \sum_m\theta_m\int_{\tau=0}^1\{f^i,f^m\}
   \vert_{x,\tau\theta}d\tau\quad
   {\rm if}\quad \t Q=\delta(\{f^i\})\,.
                    $$
   Projection of this equation on subspaces $A_{(p)}$ implies
                    $$
    f^i_{(p)}=-{\p \t Q_{(p+1)}\over\p\theta_i}+
    \hbox{terms depending on
    $f^i_{(1)},\dots f^i_{(p-1)}$}.
                \eqno ({Ap3}.8)
       $$
Hence $\delta$ is injection. Comparing this relation with
relation (Ap3.6) we see that
 the map $\delta\circ\U\colon A^{-(2)}\rightarrow A^{-(2)}$
  is bijection.
 Hence the map $\U$ is also bijection.
 For every $\{f^i\}\in\N$ the odd function
    $Q=(\delta\circ\U)^{-1}\circ$ $\delta(\{f^i\}))$
 is the unique Hamiltonian in $A^{(2)}$ required by Lemma.

\def\BarScha{1}
\def\BarSchb{2}
\def\tyutina {3}
\def\tyutinb {4}
\def \BVa {5}
\def \BVb {6}
\def\Berezin {7}
\def\Bernsteina {8}
\def\Bernsteinb {9}
\def\buttin{10}
\def \Gayduk {11}
\def\Gui  {12}
\def\JMPa {13}
\def\CMP {14}
\def\But {15}
\def\MPL{16}
\def\JMPd {17}
\def\Poin {18}
\def\leitb {19}
\def\leites{20}
\def\Shander{21}
\def\SchNucl {22}
\def\SchCMPb {23}
\def\SchCMPa{24}
\def\Vora   {25}
\def\Vor  {26}

\bigskip

\def \sm {\smallskip}

           \centerline{\bf References}

                 \medskip

 [\BarScha] M.A.Baranov, A.S.Schwarz---
  Characteristic Classes of
 Supergauge Fields.
 \noindent{\it Funkts. Analiz i ego pril.},
 {\bf 18}, No.2,  53--54, (1984).

\sm
 [\BarSchb] M.A. Baranov, A.S.Schwarz---
  Cohomologies of Supermanifolds.
 {\it Funkts. analiz i ego pril.}.
 {\bf 18}, No.3, 69--70, (1984).

\sm [\tyutina] I. A. Batalin, I. V. Tyutin---On possible
generalizations of field--antifield formalism {\it Int. J. Mod.
Phys.}, {\bf A8}, pp. 2333-2350. (1993).

\sm [\tyutinb] I. A. Batalin, I. V. Tyutin--- On the Multilevel
Field--Antifield Formalism with the Most General Lagrangian
Hypergauges. {\it Mod. Phys. Lett}, {\bf A9}, pp.1707-1712,
(1994).

\sm [\BVa] I.A.Batalin, G.A.Vilkovisky---
         Gauge algebra and Quantization.
      {\it Phys.Lett.}, {\bf 102B} pp.27--31, (1981).

\sm [\BVb] I.A.Batalin, G.A.Vilkovisky---
     Closure of the Gauge Algebra,
 Generalized Lie Equations and Feynman Rules.
  {\it Nucl.Phys.} {\bf B234}, 106--124, (1984).

\sm
 [\Berezin] F.A.Berezin
{\it Introduction to Algebra and Analysis with Anticommuting
 Va\-ri\-ables.}  Moscow, MGU (1983).
 ({\it in English}--- Introduction to Superanalysis.
  Dordrecht--Boston: D.Reidel Pub. Co., (1987)).

\sm [\Bernsteina]J.N.Bernstein, D.A. Leites---Integral forms and
the Stokes formula on supermanifolds. {\it Funkts. Analiz i ego
pril.} {\bf 11} No.1 pp. 55--56, (1977).

\sm
 [\Bernsteinb] J.N.Bernstein, D.A. Leites---
 How to integrate differential forms on
su\-per\-ma\-ni\-folds. {\it Funkts.. Analiz i ego pril.} {\bf
11} No.3 pp.70--71, (1977).

\sm
 [\buttin] C.Buttin---C.R.Acad. Sci. Paris, Ser.A--B,
  {\bf 269} A--87, (1969).

\sm
 [\Gayduk] A.V. Gayduk, O.M.Khudaverdian, A.S Schwarz---
  Integration on
  Sur\-faces in  Su\-per\-space.
  {\it Teor. Mat. Fiz.} {\bf 52}, 375--383, (1982).

\sm
 [\Gui]  V.Guillemin, S.Sternberg---{\it Geometric Asymptotics.}
AMS, Providence, Rhode Island, 1977.
 (Mathematical Surveys. Number 14)

\sm [\JMPa] O.M. Khudaverdian---
 Geometry of Superspace with Even and Odd Brackets--
    {\it J.Math.Phys}., {\bf 32}, pp. 1934--1937, (1991),
    (Preprint of Geneva University UGVA-DPT 1989/05--613).

\sm
 [\CMP] O.M.Khudaverdian---  Odd Invariant Semidensity and
 Divergence-like Operators in an Odd Symplectic Superspace,
 {\it Comm. Math.Phys}., {\bf 198}, pp.591--606, (1998).

\sm
 [\But] O.M. Khudaverdian, R.L.Mkrtchian---
   Integral Invariants of Buttin
  Bracket.
  {\it Lett. Math. Phys.}, {\bf 18}, pp. 229--234, (1989).

\sm
 [\MPL] O.M.Khudaverdian, A.P.Nersessian---
   On Ge\-omet\-ry of Ba\-ta\-lin-
   Vil\-ko\-vis\-ky
   For\-ma\-lism {\it Mod.Phys.Lett}., {\bf A8},
     No.25 pp.2377--2385, (1993).

\sm
  [\JMPd] O.M.Khudaverdian, A.P.Nersessian---
  Batalin--Vilkovisky Formalism and
  Integration Theory on Manifolds.
  {\it J. Math. Phys}., {\bf37},pp.3713-3724, (1996).

\sm
      [\Poin] O.M. Khudaverdian, A.S. Schwarz, Yu. S. Tyupkin---
Integral Invariants for Supercanonical Transformations.---
 {\it Lett. Math. Phys.}, {\bf 5}, pp. 517--522. (1981).

\sm [\leitb]. D.A.Leites--- The new Lie Superalgebras and
Mechanics.
   {\it Docl. Acad. Nauk SSSR} {\bf 236}, 804--807, (1977).

\sm
 [\leites]. D.A. Leites--- {\it The Theory of Supermanifolds.}
   Karelskij Filial AN SSSR (1983).

\sm [\Shander] V.N.Shander---
  Analogues of the Frobenius and Darboux Theorems for Supermanifolds.
  {\it Comptes rendus de l'
 Academie bulgare des Sciences}, {\bf 36}, n.3,
 309--311, (1983).

\sm [\SchNucl] A.S.Schwarz---Are the Field and Space Variables on
an equal Footing?
 {\it Nuclear Physics.}, {\bf B171},  pp.154--166, (1980).

\sm
 [\SchCMPb] A.S.Schwarz---  Supergravity, Complex Geometry and
 G-structures.
   {\it Commun. Math. Phys.}, {\bf 87}, 37--63, (1982).

\sm [\SchCMPa] A.S.Schwarz---Geometry of Batalin-Vilkovisky
Formalism.

 \noindent{\it Commun. Math. Phys.}, {\bf 155},  pp.249--260, (1993).

\sm
  [\Vora] T. Voronov, A.Zorich---Integral transformations of
 pseudodifferential forms.
 {\it Usp. Mat.Nauk } {\bf 41}, No.6, pp.167--168, (1986).

\sm
 [\Vor] T.Voronov---  Geometric Integration Theory
     on Supermanifolds.
   {\it Sov.Sci.Rev.C Math.} {\bf 9}, 1--138, (1992).

\bye